\documentclass[a4paper,11pt]{article}  
\usepackage[utf8]{inputenc} \usepackage{mathtools,amsmath,graphicx,array,amsthm,amssymb}  
\usepackage{float}  
\usepackage{subcaption} 
\usepackage[graphicx]{realboxes}  
\usepackage{hyperref} 
\usepackage{microtype}
\usepackage{bm}
\usepackage{natbib}
\usepackage{booktabs}
\usepackage{algorithm}
\usepackage{algorithmic}
\usepackage{bm}
\usepackage{ragged2e}

\newtheorem{eg}{Example}[section]  
\let\oldeg\eg  
\let\oldendeg\endeg  
\renewenvironment{eg}{%
    \oldeg  
    \normalfont}{%
    \oldendeg}  

\usepackage{authblk}

\begin{document}

\title{Multiscale Fourier Neural Operator for Inverse Wave Scattering in Highly Oscillatory Media}

\author[1]{Zilin You}
\author[1]{Zhenli Xu}
\author[2]{Wei Cai\thanks{Corresponding author, email: cai@smu.edu.}}

\affil[1]{School of Mathematical Sciences, Shanghai Jiao Tong University, Shanghai, China}
\affil[2]{Department of Mathematics, Southern Methodist University, Dallas, TX, USA}
\date{}
\maketitle

\begin{abstract}
In this paper, we propose an operator learning method based on the multiscale Fourier neural operator (MscaleFNO) for inverse medium problems of Helmholtz equations. The MscaleFNO provides a neural surrogate model with reduced spectral bias for the Helmholtz equations, mapping highly oscillatory medium profiles to scattered wavefields. A plug-and-play inversion using elucidated diffusion model is introduced to regularize the inverse solver based on least squares of data misfits. Numerical results for partial aperture inversion of oscillatory two-dimensional media demonstrate the advantage and effectiveness of MscaleFNO for accurate reconstruction of highly oscillatory medium properties.  
\end{abstract}
\begin{quote}
\end{quote}

\section{Introduction}
Inverse medium problems arise in a wide range of applications, including radar imaging, geophysical exploration, and medical diagnostics, where the objective is to reconstruct unknown material parameters from observed wavefield measurements. Among existing approaches, the full waveform inversion (FWI) has emerged as a powerful framework that exploits full information content of wavefields to recover subsurface or interior structures \cite{Tarantola1984, Virieux2009}. Unlike classical travel-time tomography, which relies only on the kinematics of wave propagation, FWI is formulated as a PDE-constrained optimization problem that minimizes the mismatch between observed and simulated data, enabling higher-resolution reconstructions of fine-scale features \cite{Virieux2009}. Despite its success, FWI faces fundamental challenges at high frequencies, stemming from the non-convexity of the misfit functional and the prohibitive computational cost of repeated PDE solves.

The rapid development of deep learning has introduced data-driven alternatives to classical iterative approaches for inverse problems. In the context of FWI, early works focused on end-to-end supervised mappings from wavefield data to medium parameters using convolutional architectures, including InversionNet \cite{Wu2020} and VelocityGAN \cite{Zhang2019}, as well as hybrid approaches such as UPFWI \cite{Jin2022} that incorporate differentiable forward models. While these methods offer significant acceleration at inference, they often lack physical consistency, generalize poorly to out-of-distribution media, and fail to recover high-frequency structures due to spectral bias \cite{Xu_2025_overview}. Benchmark studies such as OpenFWI \cite{Deng2022} and large-scale models \cite{BigFWI2024} further highlight the persistent trade-off between computational efficiency and reconstruction fidelity. 

Operator learning provides a principled alternative by approximating solution operators of parametric PDEs. Representative architectures include Fourier Neural Operators (FNO) \cite{Li2021}, DeepONet \cite{chen1995,Lu2021}, and physics-informed variants such as PINO \cite{PINO2024}, along with extensions for complex geometries and multiscale representations \cite{GNO2020,WNO2023,PINTO2025}, and particle-system-inspired formulations 
that offer improved generalizability for PDE approximation \cite{Ji2025}. These approaches have demonstrated strong performance in learning forward solution operators for PDEs \cite{Kovachki2023}, and have been increasingly employed in inverse problems either as differentiable surrogate forward models within optimization frameworks \cite{Cao2025DiffANO,Liu2025Neumann}, or as direct inverse mappings from observations to parameters \cite{Khoo2019,Molinaro2023,Long2024,Rahman2022}. However, the latter often suffers from poor generalization, heavy reliance on paired training data, and loss of fine-scale details due to the ill-posedness of the inverse problem. A fundamental limitation shared by many neural operators is their bias toward low-frequency components, which hampers their ability to represent highly oscillatory wavefields \cite{Xu2025,Rahaman2019}. Recent efforts to address this include hierarchical and diffusion-enhanced operator architectures \cite{HANO2024,oommen2024}. More recently, the multiscale Fourier neural operator (MscaleFNO) \cite{You2026} has been proposed to overcome this limitation through scale-adaptive representations, achieving significantly improved approximation of high-frequency solutions. Nevertheless, deploying such operators within a statistically principled inversion framework that accounts for the ill-posedness of the problem remains an open challenge.

Diffusion models have emerged as powerful generative frameworks for learning complex high-dimensional data distributions. Building on DDPM \cite{Ho2020}, subsequent developments including score-based models \cite{Song2021}, DDIM \cite{Song2020b}, EDM \cite{Karras2022}, and latent diffusion \cite{Rombach2022} have established diffusion models as state-of-the-art tools for image generation and flexible prior modeling. These capabilities have been increasingly leveraged in inverse problems, where diffusion models serve as learned priors for regularizing ill-posed reconstructions. Representative approaches include posterior sampling methods such as DPS \cite{Chung2023DPS}, DDRM \cite{Kawar2022}, and DiffPIR \cite{Zhu2023}. A complementary line of work integrates diffusion-based denoisers into iterative optimization schemes, including Plug-and-Play (PnP) \cite{Venkatakrishnan2013} and Regularization by Denoising (RED) \cite{Romano2017} frameworks. Recent extensions such as RED-DiffEq \cite{REDDiffEq2026} have demonstrated strong performance in PDE-constrained inverse problems including FWI. However, all such approaches still require repeated evaluations of the forward PDE solver to enforce data consistency. In high-frequency regimes governed by the Helmholtz equation, where extremely fine discretization is necessary, such repeated solves become computationally prohibitive, making the forward solver the dominant bottleneck and fundamentally limiting the scalability of diffusion-regularized inversion for high-resolution FWI.

The core challenge in high-frequency inverse scattering is fundamentally computational. Resolving fine-scale structures demands large wavenumbers, which render traditional numerical solvers prohibitively expensive for the repeated PDE evaluations required in iterative optimization. Neural operator surrogates offer a potential escape from this expense, but standard architectures such as FNO suffer from spectral bias, failing to accurately represent highly oscillatory wavefields. Together, these limitations leave high-frequency FWI without an efficient and physically consistent solution, while the ill-posedness of the inverse problem remains inadequately addressed.
To address this challenge, we propose a full waveform inversion framework built upon the MscaleFNO \cite{You2026}. By leveraging scale-adaptive representations, MscaleFNO overcomes the spectral bias of standard FNOs and enables accurate wavefield prediction in high-wavenumber regimes, serving as an efficient differentiable surrogate that eliminates the need for repeated expensive PDE solves. To further regularize the ill-posed inversion, we incorporate a pre-trained elucidated diffusion model (EDM) as a plug-and-play prior, injecting learned structural priors into the iterative reconstruction without additional forward solver evaluations. The resulting framework achieves high-frequency FWI that is simultaneously efficient, accurate, and well-regularized.

The main contributions of this work are as follows:
\begin{itemize}
\item We demonstrate that MscaleFNO serves as a superior surrogate solver for high-wavenumber Helmholtz equations, achieving faster inference, improved numerical stability, and higher solution accuracy compared to standard FNO architectures, owing to its scale-adaptive spectral representations.

\item We embed MscaleFNO as a differentiable surrogate forward model within a FWI framework for high-frequency Helmholtz inverse medium problems, replacing repeated expensive PDE solves and enabling efficient gradient-based optimization.

\item We validate the proposed framework under noisy and limited-aperture measurement settings, demonstrating that the integration of an EDM-based plug-and-play diffusion prior further improves reconstruction accuracy, robustness to noise, and recovery of fine-scale structural details compared to baseline methods.
\end{itemize}

The remainder of this paper is organized as follows. Section 2 briefly reviews the formulation of full waveform inversion. Section 3 presents the proposed framework, detailing the MscaleFNO-based surrogate forward model and the integration of the EDM-based plug-and-play diffusion prior. Section 4 presents numerical experiments validating the surrogate forward model, evaluating reconstruction performance under noisy and limited-aperture settings, and conducting ablation studies. Finally, Section 5 concludes the paper.

\section{Full waveform inversion}
The FWI is typically implemented using gradient-based methods such as the adjoint-state approach \cite{Plessix2006}, which require repeated solutions of the forward wave equation. In high-frequency regimes, the wavefield becomes highly oscillatory, necessitating mesh resolutions proportional to the wavelength. This leads to large, indefinite, and poorly conditioned linear systems \cite{Amundsen2005, Dorn2000}, significantly increasing computational cost. Moreover, the achievable resolution is fundamentally limited by diffraction, with features below approximately half a wavelength being difficult to resolve \cite{Nobel2014, Courjon1994}. Achieving higher resolution thus requires information of higher frequencies, which further exacerbates both computational burden and ill-posedness \cite{Hahner2001}. In addition, numerical discretizations suffer from the pollution effect in the high-frequency regime, requiring prohibitively fine meshes \cite{Ihlenburg1997}.

The propagation of acoustic waves in inhomogeneous media is governed by the Helmholtz equation
\begin{equation}
\Delta u(\bm x) + k^2 \big( 1 + \alpha(\bm x) \big) u(\bm x) = f(\bm x), 
\qquad \bm x \in \Omega,
\label{helm_eq}
\end{equation}
where $k$ denotes the wavenumber in a homogeneous background medium, 
$\alpha(\bm x)$ represents the relative perturbation of the medium, and $f(\bm x)$ is the source term.
The \textbf{inverse medium problem} aims to reconstruct $\alpha(\bm x)$ from observations of the wavefield $u(\bm x)$ \cite{Colton2000,Nedelec2001}. In the FWI setting, the total wavefield is matched directly to observed data. The total field can be decomposed as
\begin{equation}
u(\bm x) = u^{inc}(\bm x) + u^s(\bm x),
\end{equation}
where $u^{inc}(\bm x)$ is the incident field and $u^s(\bm x)$ is the scattered field.

We consider point sources of the form
$ 
f(\bm x) = \delta(\bm x - \bm x_s)
$
located at position $\bm x_s$. In a homogeneous medium, the incident field satisfies
\begin{equation}
\Delta u^{inc}(\bm x) + k^2 u^{inc}(\bm x) = \delta(\bm x - \bm x_s),
\end{equation}
which corresponds to a spherical (or cylindrical) wave.
The scattered field $u^s(\bm x)$ satisfies
\begin{equation}
\Delta u^s(\bm x) + k^2 (1 + \alpha(\bm x)) u^s(\bm x) = -k^2 \alpha(\bm x) u^{inc}(\bm x),
\end{equation}
together with the Sommerfeld radiation condition
\begin{equation}
\lim_{|\bm x|\to\infty} |\bm x|^{\frac{d-1}{2}}
\left(
\frac{\partial u^s}{\partial |\bm x|} - i k u^s
\right) = 0,
\end{equation}
which ensures that the scattered field is outgoing.

To probe the medium, we consider multiple sources and receivers. Let
$
\{\bm x_s^m\}_{m=1}^M$ and
$
\{\bm x_r^n\}_{n=1}^N$ in $\Omega$
denote the source and receiver locations. For each source $\bm x_s^m$, we define
$f_m = \delta(\bm x - \bm x_s^m)$ and denote the corresponding wavefield by
$u(\bm x; \bm x_s^m)$.
We define the forward Helmholtz operator
$
\mathcal{K} : (\alpha, f) \mapsto u,
$
which maps the medium perturbation and source to the corresponding wavefield.
The noisy measurements are collected in a complex-valued observation matrix
$
\mathbf{Y}^\delta = 
\{ y_{m,n}^\delta \}_{m=1,n=1}^{M,N},
$
with entries
\begin{equation}
y_{m,n}^\delta 
= 
\mathcal{K}(\alpha, f_m)(\bm x_r^n)
+ \epsilon_{m,n},
\end{equation}
where $\epsilon_{m,n}$ denotes complex-valued additive noise.

The inverse problem is to reconstruct $\alpha(\bm x)$ from the noisy multi-source observations $\mathbf{y}^\delta$. Due to the ill-posedness of the problem, regularization is required to ensure stability. A standard formulation is the PDE-constrained optimization problem
\begin{equation}
\min_{\alpha(\bm x)} 
\; \mathcal{J}(\alpha)
=
\frac{1}{MN}
\sum_{m=1}^{M}\sum_{n=1}^{N}
\left\|
\mathcal{K}(\alpha, f_m)(\bm x_r^n) - y_{m,n}^\delta
\right\|_2^2
+ \mathcal{R}(\alpha),
\end{equation}
where $\mathcal{R}(\alpha)$ is a regularization term encoding prior knowledge on the medium structure. In this work, this term is realized through a learned plug-and-play diffusion prior, as detailed in the next section. 

\section{Operator learning for full waveform inversion}  
\subsection{The algorithm framework}

We consider the inverse medium problem of recovering the perturbation field $\alpha(\bm x)$ from noisy wavefield observations $\bm Y^\delta=\{y_{m,n}^{\delta}\}$, where $m=1,\dots,M$ indexes the source and $n=1,\dots,N$ indexes the receiver. A conventional regularized formulation can be written as
\[
\min_{\alpha} \ \mathcal{L}(\alpha)+\mathcal{R}(\alpha),
\]
where \(\mathcal{L}(\alpha)\) measures the discrepancy between the predicted and observed wavefields, and \(\mathcal{R}(\alpha)\) denotes a regularization term. In this work, rather than explicitly specifying and differentiating such a regularization functional, we employ a learned plug-and-play diffusion prior to periodically guide the inversion iterates.

Let $\mathcal{K}$ denote the forward Helmholtz operator that maps the medium perturbation $\alpha(\bm x)$ and a source term $f_m$ to the corresponding wavefield. In conventional PDE-constrained inversion, repeated evaluations of the forward Helmholtz solver and its adjoint are required during the optimization process, which becomes computationally expensive, especially in high-frequency regimes. To alleviate this cost, we replace the expensive forward operator $\mathcal{K}$ with a pretrained neural operator surrogate $\mathcal{K}_\theta$. During inversion, the parameters $\theta$ of the neural operator are fixed, and only the medium perturbation $\alpha$ is optimized. Since $\mathcal{K}_\theta$ is implemented as a differentiable neural network, gradients with respect to $\alpha$ can be efficiently computed by automatic differentiation.

The inversion is initialized from a Gaussian random field,
\begin{equation}
\alpha^{(0)}
\sim
\mathcal{N}(0,\sigma_{\max}^{2}I),
\end{equation}
which is treated as an optimizable tensor with gradient tracking enabled. At the $i$th inversion iteration, the current estimate $\alpha^{(i)}$ is fed into the neural operator together with all source terms in a batched manner. Specifically, we write the batched neural-operator prediction as
\begin{equation}
\widehat{\bm Y}^{(i)}
=
\mathcal{K}_{\theta}^{\mathrm{batch}}
\left(
\{(\alpha^{(i)}, f_m)\}_{m=1}^{M}
\right):=
\{\mathcal{K}_{\theta}
(\alpha^{(i)}, f_m)
\}_{m=1}^{M},
\end{equation}
where the predicted observation corresponding to the $m$th source and the $n$th receiver is given by
\begin{equation}
\widehat{Y}_{m,n}^{(i)}
=
\mathcal{K}_{\theta}
\left(
\alpha^{(i)}, f_m
\right)(\bm x_r^n).
\end{equation}
In practice, the current medium estimate \(\alpha^{(i)}\) is broadcast to match the batch dimension of the source terms, so that the medium-source pairs \(\{(\alpha^{(i)}, f_m)\}_{m=1}^{M}\) can be processed simultaneously by the neural operator.  This batched evaluation avoids source-by-source forward simulations and enables efficient parallel computation on modern accelerators.

The data fidelity term is then defined over the entire source-receiver observation set as
\begin{equation}
\mathcal{L}(\alpha^{(i)})
=
\frac{1}{MN}
\sum_{m=1}^{M}
\sum_{n=1}^{N}
\left\|
\widehat{Y}_{m,n}^{(i)} - y_{m,n}^{\delta}
\right\|_2^2.
\end{equation}
Instead of explicitly deriving the adjoint-state gradient, we compute the gradient of the data fidelity term with respect to the medium variable by backpropagating through the differentiable surrogate model:
\begin{equation}
g^{(i)}
=
\nabla_{\alpha}
\mathcal{L}(\alpha^{(i)}).
\end{equation}
This gradient is obtained through the computational graph
\begin{equation}
\alpha^{(i)}
\longrightarrow
\mathcal{K}_{\theta}^{\mathrm{batch}}
\longrightarrow
\widehat{\bm Y}^{(i)}
\longrightarrow
\mathcal{L}(\alpha^{(i)}),
\end{equation}
Therefore, the proposed formulation avoids the need for explicit Jacobian construction or adjoint Helmholtz solvers.
The medium perturbation $\alpha^{(i)}$ is then updated by an Adam optimizer using the automatically computed gradient. We denote this data-consistency update by
\begin{equation}
\tilde{\alpha}^{(i)}
=
\mathrm{AdamStep}
\left(
\alpha^{(i)}, g^{(i)}; \eta
\right),
\end{equation}
where $\eta$ is the learning rate. This update reduces the mismatch between the predicted observations and the measured data, thereby enforcing consistency with the physical measurements.

Due to the ill-posedness of the inverse problem, data consistency alone may lead to unstable reconstructions or nonphysical artifacts. Instead of prescribing an explicit handcrafted regularization functional $\mathcal{R}(\alpha)$, we introduce a learned prior based on diffusion models to implicitly characterize the distribution of plausible medium perturbations. In particular, we employ an Elucidated Diffusion Model as a learned plug-and-play diffusion prior. The diffusion model is not applied after every gradient update. Rather, it is invoked periodically to correct the current reconstruction and guide it toward the learned medium prior. The update rule is written as
\begin{equation}
\alpha^{(i+1)}
=
\begin{cases}
\mathrm{EDM\_sampler}
\left(
\tilde{\alpha}^{(i)}, \{\sigma_s\}_{s=0}^{S}
\right),
&
\text{if } (i+1)\bmod T = 0, \\[4pt]
\tilde{\alpha}^{(i)},
&
\text{otherwise},
\end{cases}
\end{equation}
where $T$ denotes the diffusion update period, $S$ is the number of denoising steps used by the EDM sampler, and $\{\sigma_s\}_{s=1}^{S}$ denotes the noise schedule. This periodic prior update allows several data-consistency steps to sufficiently reduce the observation mismatch before applying the diffusion-based correction. Compared with applying the diffusion model at every iteration, the periodic strategy reduces computational overhead and mitigates excessive regularization from the learned prior.

The overall framework therefore alternates between efficient data-consistency optimization and periodic diffusion-prior correction. The data consistency step enforces agreement with the measured wavefields through the differentiable neural operator, while the diffusion step suppresses artifacts and constrains the reconstruction to the learned distribution of physically plausible media. The proposed framework avoids repeated numerical PDE solves during inversion by combining a pretrained neural operator surrogate with automatic-differentiation-based gradient computation. Batched multi-source evaluation and periodic EDM correction further improve computational efficiency and reconstruction quality. The overall procedure is summarized in Algorithm~\ref{alg:pnp_diffusion}.

\begin{algorithm}[H]
\caption{Diffusion-Prior-Guided Inversion with Neural Operator}
\label{alg:pnp_diffusion}
\begin{algorithmic}[1]

\STATE \textbf{Input:} Observed data $\bm Y^\delta=\{y_{m,n}^\delta\}$, sources $\{f_m\}_{m=1}^{M}$, receiver locations $\{\bm x_r^n\}_{n=1}^{N}$,
\STATE \hspace{\algorithmicindent} pretrained neural operator $\mathcal{K}_\theta$, EDM sampler, learning rate $\eta$, total iterations $I$, diffusion update period $T$, EDM noise schedule $\{\sigma_s\}_{s=0}^{S}$
\STATE \textbf{Initialize:} Sample $\alpha^{(0)} \sim \mathcal{N}(0,\sigma_{\max}^2 I)$ and set $\alpha^{(0)}$ as an optimizable tensor
\STATE Freeze the neural operator parameters $\theta$ and initialize the Adam optimizer for $\alpha$

\FOR{$i = 0,1,\dots,I-1$}

    \STATE \textbf{Batched forward prediction:}
    \begin{equation*}
    \widehat{\bm Y}^{(i)}
    =
    \mathcal{K}_{\theta}^{\mathrm{batch}}
    \left(
    \alpha^{(i)}, \{f_m\}_{m=1}^{M}
    \right)
    \end{equation*}

    \STATE \textbf{Compute data fidelity loss:}
    \begin{equation*}
    \mathcal{L}(\alpha^{(i)})
    =
    \frac{1}{MN}
    \left\|
    \widehat{\bm Y}^{(i)}
    -
    \bm Y^\delta
    \right\|_F^2
    \end{equation*}

    \STATE \textbf{Automatic differentiation:}
    \begin{equation*}
    g^{(i)}
    =
    \nabla_{\alpha}
    \mathcal{L}(\alpha^{(i)})
    \end{equation*}

    \STATE \textbf{Data-consistency update using Adam:}
    \begin{equation*}
    \tilde{\alpha}^{(i)}
    =
    \mathrm{AdamStep}
    \left(
    \alpha^{(i)}, g^{(i)}; \eta
    \right)
    \end{equation*}

    \IF{$(i+1)\bmod T = 0$}
        \STATE \textbf{Periodic prior update using EDM:}
        \begin{equation*}
        \alpha^{(i+1)}
        =
        \mathrm{EDM\_sampler}
        \left(
        \tilde{\alpha}^{(i)}, \{\sigma_s\}_{s=0}^{S}
        \right)
        \end{equation*}
    \ELSE
        \STATE \textbf{No diffusion update:}
        \begin{equation*}
        \alpha^{(i+1)}
        =
        \tilde{\alpha}^{(i)}
        \end{equation*}
    \ENDIF

\ENDFOR

\STATE \textbf{Output:} Reconstructed medium perturbation $\alpha^{(I)}$

\end{algorithmic}
\end{algorithm}

\subsection{Neural operator forward surrogate}

We learn a surrogate forward operator \(\mathcal{K}_{\theta}\) that predicts the total wavefield \(u(x)\) from the medium perturbation \(\alpha(x)\) and the source term \(f(x)\), with spatial coordinates concatenated as positional inputs in the neural representation:
\[
\mathcal{K}_{\theta} : (\alpha, f) \mapsto u.
\]
According to the diffraction limit, accurate reconstruction of high-frequency structures requires sufficiently large wavenumbers $k$, so that the wave can adequately probe fine-scale features of the medium. However, increasing $k$ leads to highly oscillatory solutions, which poses a fundamental challenge for neural operator learning due to spectral bias. Specifically, conventional operator networks (e.g., FNO and DeepONet) tend to favor low-frequency components and struggle to accurately represent high-frequency oscillations.

\noindent\textbf{Fourier Neural Operator.}
The FNO learns mappings between function spaces by parameterizing integral operators in the spectral domain. Given a generic input function $a(\bm x)$ defined on domain $D \subset \mathbb{R}^d$, the FNO constructs an operator $G_\theta$ such that
\begin{equation}
u(\bm x) = G_\theta(a(\bm x),\bm x),
\end{equation}
where the coordinate $\bm x$ is concatenated with $a(\bm x)$ as an additional input channel, equivalently serving as an identity coordinate map on $D$.

\begin{figure}[t!]
    \centering
    \includegraphics[width=0.85\textwidth]{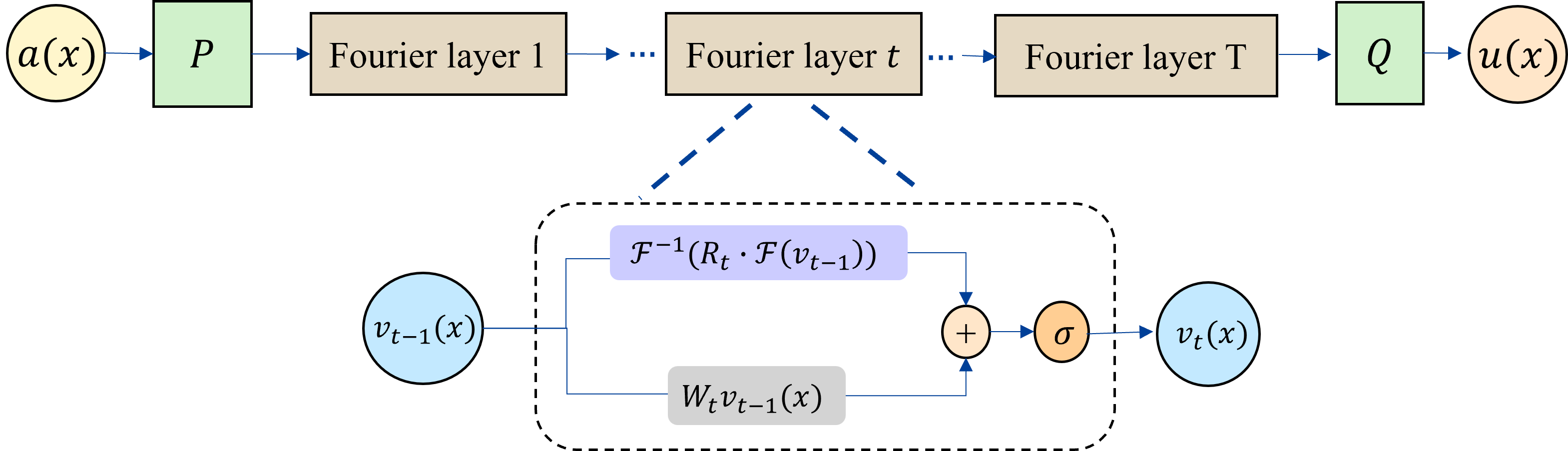}
    \caption{
        Architecture of the FNO network. The model consists of a lifting layer, multiple Fourier layers, and a projection layer, where each Fourier layer combines local transformations in the spatial domain with global convolutions in the spectral domain.
    }
    \label{fig:fno_architecture}
\end{figure}

The architecture consists of three components: a lifting operator $P$, a sequence of Fourier layers, and a projection operator $Q$. The input is first lifted to a higher-dimensional feature space:
\begin{equation}
v_0(\bm x) = P\big(a(\bm x), \bm x\big).
\end{equation}
Each Fourier layer updates the feature representation via
\begin{equation}
v_{t+1}(\bm x) = \sigma\Big( W_t v_t(\bm x) + \mathcal{F}^{-1}\big(R_t \cdot \mathcal{F}(v_t)\big)(\bm x) \Big),
\end{equation}
where $W_t$ is a learnable linear transformation acting pointwise in the spatial domain, and $R_t$ is a complex-valued weight tensor acting on Fourier modes. The second term corresponds to a global convolution operator implemented efficiently in the spectral domain. This formulation combines local interactions (via $W_t$) and global dependencies (via Fourier convolution).
After $T$ layers, the output is obtained through a projection operator:
\begin{equation}
u(\bm x) = Q\big(v_T(\bm x)\big).
\end{equation}
To improve computational efficiency, the Fourier representation is truncated by retaining only the lowest $k_{\max}$ modes, i.e., $R_t$ acts only on a subset of Fourier coefficients while higher-frequency modes are discarded. Although this truncation reduces complexity, it introduces a bias toward smooth, low-frequency functions, limiting the ability of standard FNOs to represent highly oscillatory solutions. The total architecture of FNO is illustrated in Figure~\ref{fig:fno_architecture}.
As a result, standard FNOs may struggle to accurately represent highly oscillatory solutions, particularly in high-frequency regimes where important information is carried by high-wavenumber components. This motivates the use of multiscale operator learning.

\noindent\textbf{Multiscale Fourier Neural Operator.} To overcome this limitation, we employ an MscaleFNO~\cite{You2026} to approximate $\mathcal{K}_\theta$. The model consists of multiple Fourier operator branches at different scales, which helps represent both low-frequency trends and high-frequency oscillatory components. The output of the operator $\mathcal{K}_\theta$ is defined as
\begin{equation}
    u(\bm x) = \sum_{i=1}^{N}
      \gamma_i\,\mathrm{FNO}_i\big([\,c_i\alpha(\bm x),\,c_i f(\bm x),\,c_i \bm x\,]\big),
\end{equation}
where $\mathrm{FNO}_i$ denotes the Fourier neural operator at scale $i$, 
$c_i$ is the scaling factor, and $\gamma_i$ represents learnable weights. Each branch $\mathrm{FNO}_i$ adopts the FNO architecture described above.

\begin{figure}[t!]
    \centering
    \includegraphics[width=0.9\textwidth]{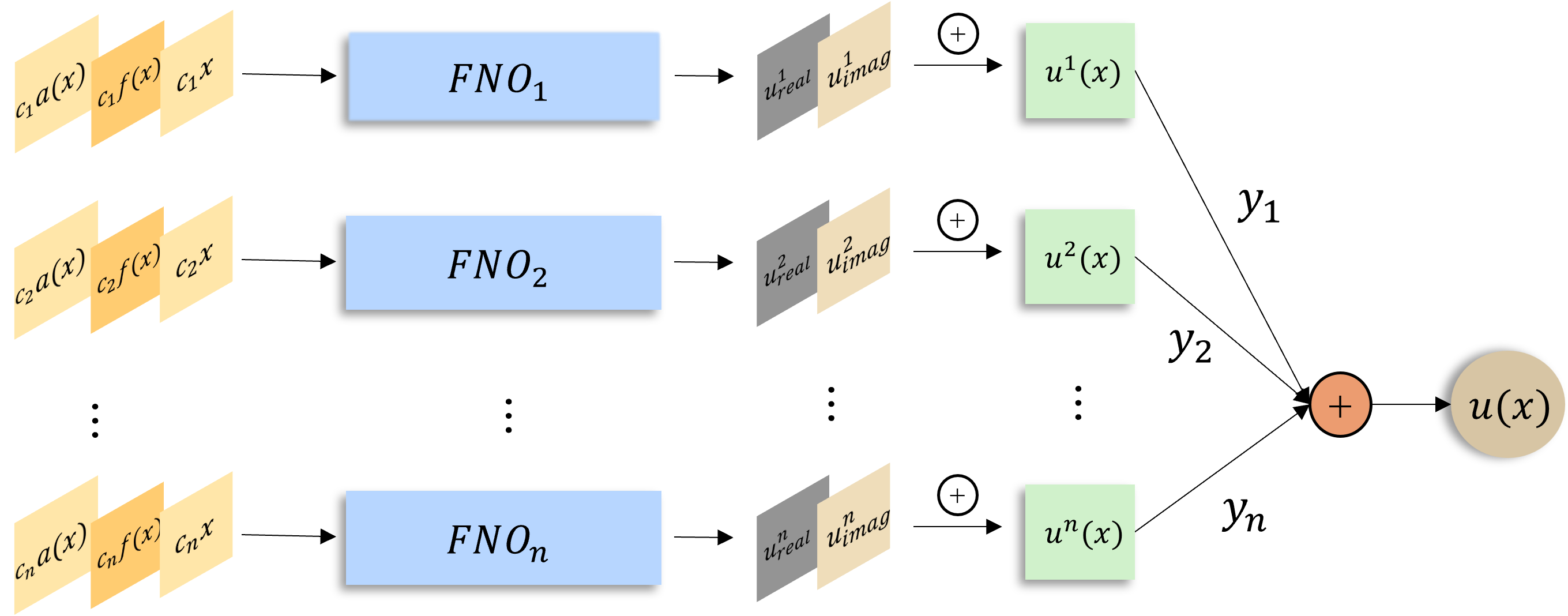}
    \caption{
        MscaleFNO architecture for approximating the operator $\mathcal{K}_\theta$.
    }
    \label{fig:mscalefno_architecture}
\end{figure}

In practice, the wave field $u(\bm x)$ is complex-valued. We represent its real and imaginary parts as two separate channels. The network therefore outputs two channels corresponding to \(\Re(u)\) and \(\Im(u)\), which are combined as
\[
{u}(\bm x) = {u}_{real}(\bm x) + i\,{u}_{imag}(\bm x).
\]
The architecture of MscaleFNO is illustrated in Figure~\ref{fig:mscalefno_architecture}.

\begin{figure}[t!]
\centering
\begin{minipage}{0.46\textwidth}
\centering
\includegraphics[width=\linewidth]{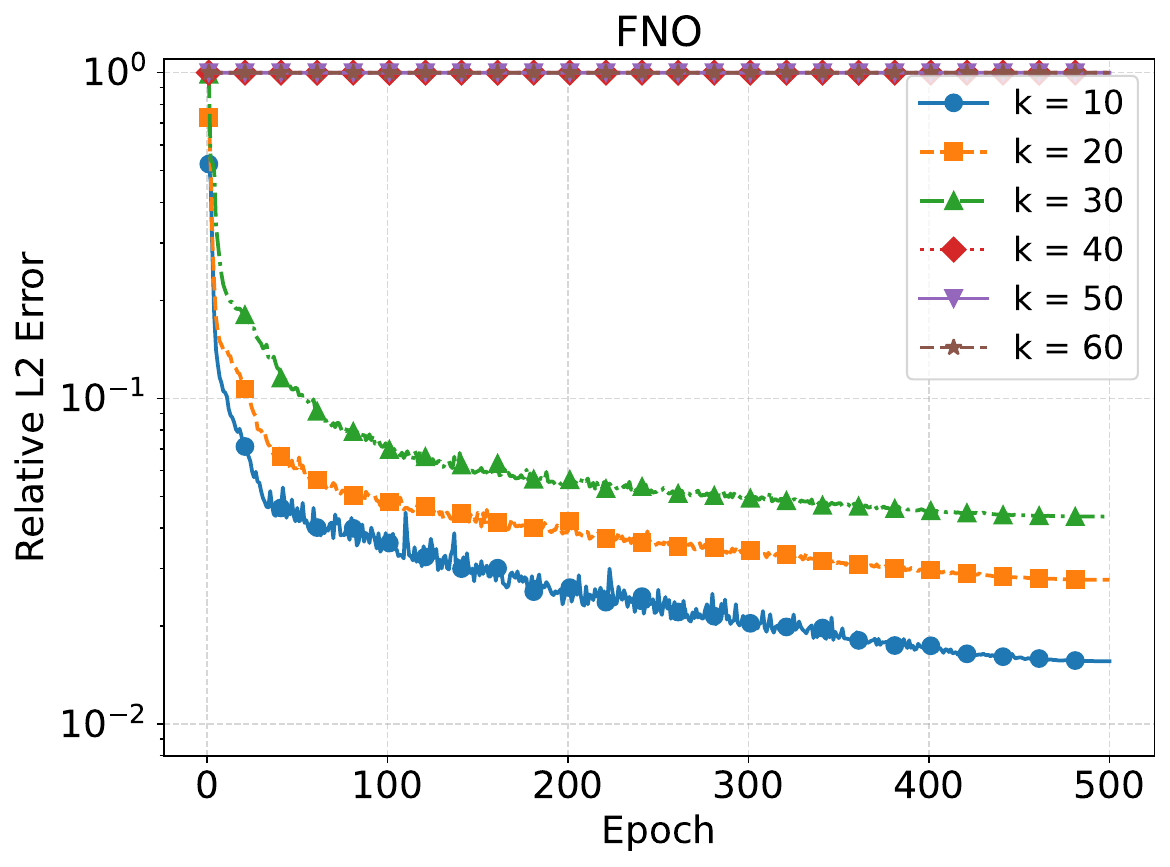}
\\ (a) FNO
\end{minipage}
\hspace{0.05\textwidth}
\begin{minipage}{0.46\textwidth}
\centering
\includegraphics[width=\linewidth]{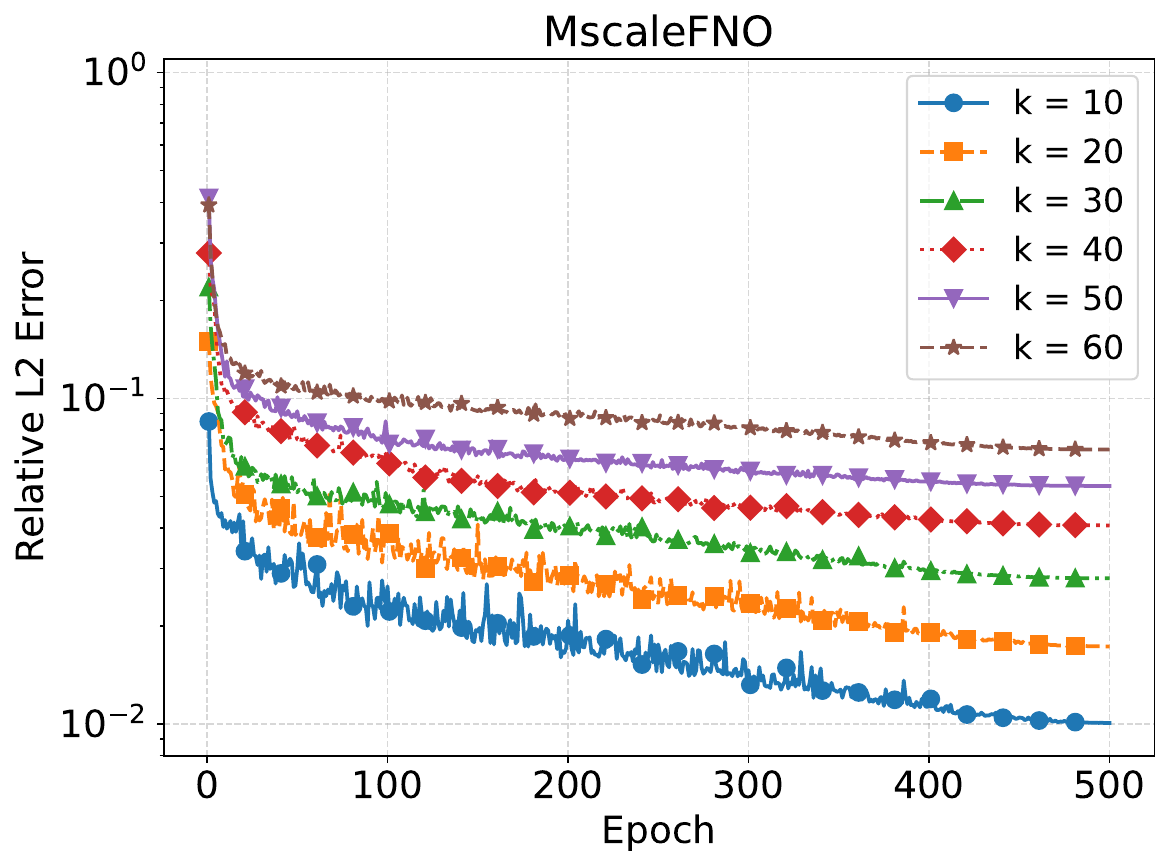}
\\ (b) MscaleFNO
\end{minipage}
\caption{Training error of FNO and MscaleFNO under different wavenumbers $k$.}
\label{fig:forward_loss_2d}
\end{figure}

Figure~\ref{fig:forward_loss_2d} shows the error convergence curves of FNO and MscaleFNO when learning the Helmholtz forward operator at different wavenumbers, where the experimental configuration follows that of 
Section~\ref{sec:forward}. For small $k$, FNO achieves reasonable accuracy with fast convergence. However, as $k$ increases, FNO exhibits noticeably slower convergence and larger residual errors, indicating difficulty in capturing high-frequency components. In contrast to FNO, MscaleFNO maintains stable convergence and consistently achieves lower training errors across all tested frequencies, demonstrating its robustness in the high-frequency regime.


\noindent\textbf{Training Data.}
To train the model, we generate a dataset of physically consistent samples
\begin{equation}
\{\,\alpha_r(\bm x),\,f_m(\bm x),\,u_{m}^r(\bm x)\,\}_{m,r=1}^{M,R},
\end{equation}
where $u_{m}^r(\bm x)$ denotes the total wave field corresponding to medium $\alpha_r$ and source $f_m$. The wavefields 
$u_{m}^r(\bm x)$ are computed by numerically solving the Helmholtz equation~\eqref{helm_eq} with a finite element solver for each medium-source pair.

\noindent \textbf{Loss Function and Training.}
The network is trained using a complex-valued relative $L_2$ loss:
\begin{equation}
    \mathcal{L}
    = \frac{1}{MR}\sum_{m,r}
      \frac{\left\|
        \mathcal{K}_\theta(\alpha_r,f_m)(\bm x) -u_{m}^r(\bm x)
      \right\|_{2,\mathbb{C}}^{2}}
      {\left\| u_{m}^r(\bm x)\right\|_{2,\mathbb{C}}^{2}},
\end{equation}
where $\|\cdot\|_{2,\mathbb{C}}$ denotes the complex $L_2$ norm.

All model parameters, including the scale weights $\gamma_i$ and scale factors $c_i$ are trainable. After training, $\mathcal{K}_\theta$ serves as an efficient surrogate for the forward wave propagation operator, enabling fast and accurate predictions across different media configurations and source conditions.

\subsection{Plug-and-play diffusion prior}
To model the prior distribution of the medium parameter, we adopt the EDM~\cite{Karras2022}, which consists of a noise-conditional denoiser and a multi-step sampling procedure.

\noindent\textbf{EDM Denoiser.}
The core component of EDM is a noise-conditional denoiser $D_\theta(\bm z,\sigma)$, which takes a noisy input $\bm z$ at noise level $\sigma$ and predicts the corresponding clean signal:
\begin{equation}
\hat{\bm z}_0 = D_\theta(\bm z, \sigma).
\end{equation}
Following the EDM parameterization, the denoiser is defined as
\begin{equation}
D_\theta(\bm z,\sigma)
=
c_{\text{skip}}(\sigma)\,\bm z
+
c_{\text{out}}(\sigma)\,
F_\theta\big(c_{\text{in}}(\sigma)\,\bm z,\sigma\big),
\end{equation}
where $F_\theta$ is a neural network and $c_{\text{in}}, c_{\text{skip}}, c_{\text{out}}$ are scaling functions depending on $\sigma$.

\noindent\textbf{Training of the Denoiser.}
The denoiser is trained using noise perturbation. Given clean data $\bm z_0$, noisy samples are generated as
\begin{equation}
\bm z = \bm z_0 + \sigma \epsilon, \quad \epsilon \sim \mathcal{N}(0, I),
\end{equation}
where the noise level $\sigma$ is sampled from a prescribed distribution. The training objective is
\begin{equation}
\mathcal{L}_{\mathrm{EDM}}
=
\mathbb{E}_{\bm z_0,\epsilon,\sigma}
\left[
w(\sigma)\,\|D_\theta(\bm z,\sigma) - \bm z_0\|_2^2
\right],
\end{equation}
where $w(\sigma)$ balances contributions across different noise levels.

\noindent\textbf{EDM Sampler.}
While the denoiser performs single-step reconstruction, EDM defines a sampler that progressively refines a sample across multiple noise levels. The update direction is given by
\begin{equation}
d(\bm z,\sigma) = \frac{\bm z - D_\theta(\bm z,\sigma)}{\sigma},
\end{equation}
which approximates the score function.
A decreasing noise schedule is defined as
\begin{equation}
\sigma_{\max} = \sigma_S > \sigma_{S-1} > \cdots > \sigma_0 = \sigma_{\min}.
\end{equation}
Starting from a noisy input $\bm z$, one iteratively updates the sample using an Euler step:
\begin{equation}
\bm z_{s} = \bm z_{s+1} + (\sigma_{s} - \sigma_{s+1})\, d(\bm z_{s+1},\sigma_{s+1}),
\end{equation}
optionally followed by a second-order (Heun) correction.
Compared to a single denoising step, the multi-step sampling process helps large noise levels capture global structures, while small noise levels recover fine details. This coarse-to-fine refinement improves reconstruction quality.

\noindent\textbf{EDM Prior Operator.}
At inference, EDM performs multi-step denoising along a decreasing noise schedule $\{\sigma_s\}$, which defines a prior operator
\begin{equation}
\bm z_0 = \mathrm{EDM\_sampler}(\bm z, \{\sigma_s\}),
\end{equation}
which maps a noisy input $\bm z$ to a refined sample through iterative denoising. This operator corresponds to the periodic correction step in Algorithm~\ref{alg:pnp_diffusion}, where it is applied every $T$ iterations to guide the medium estimate toward the learned prior distribution.

\section{Numerical results}
\label{sec:experiment setup}
\begin{figure}[t!]
\centering
\includegraphics[width=0.6\linewidth]{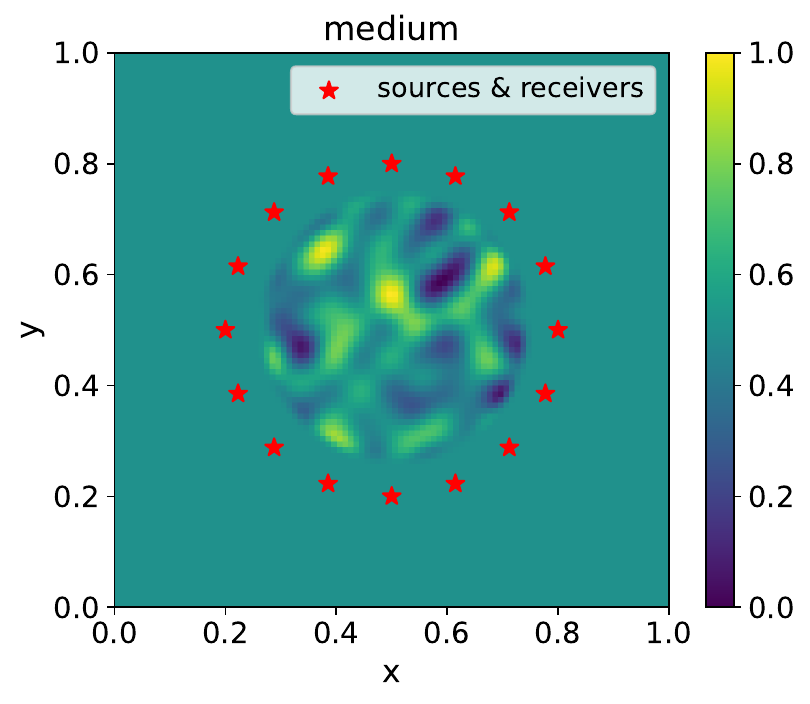}
\caption{Representative medium sample $\alpha(x,y)$ with source and receiver locations.}
\label{fig:dataset_2d_medium}
\end{figure}

In this section, we consider the solution of two-dimensional Helmholtz equations to evaluate the performance of the proposed method under various configurations. Specifically, we solve
\begin{equation}
\Delta u(x,y) + k^2 \big(1 + \alpha(x,y)\big) u(x,y) = f(x,y),
\quad (x,y)\in[0,1]^2,
\label{helm_2d_eq}
\end{equation}
where the medium perturbation is supported within a disk centered at $(0.5,0.5)$ with radius $0.25$, and vanishes outside:
\[
\alpha(x,y)=0, \quad \|(x,y)-(0.5,0.5)\| > 0.25.
\]
Inside the disk, $\alpha(x,y)$ is defined via a truncated Fourier expansion:
\begin{equation}
\begin{aligned}
\alpha(x,y)=\xi\sum_{m,n=0}^{10} \big[
& a_{mn}\sin(2\pi m x)\sin(2\pi n y) + b_{mn}\sin(2\pi m x)\cos(2\pi n y)\\
&+ c_{mn}\cos(2\pi m x)\sin(2\pi n y) + d_{mn}\cos(2\pi m x)\cos(2\pi n y)
\big],
\end{aligned}
\end{equation}
where $a_{mn}, b_{mn}, c_{mn}, d_{mn} \sim \mathrm{Uniform}(-1,1)$ and $\xi = 0.02$.
The source locations are uniformly distributed along a circle centered at $(0.5,0.5)$ with radius $0.3$, i.e.,
$\|(x,y)-(0.5,0.5)\|=0.3.$ 
The receivers are placed at grid points closest to the same circle used for the source distribution. Figure~\ref{fig:dataset_2d_medium} shows a representative example of the medium along with the corresponding source and receiver configurations.

\begin{figure}[t!]
\centering
\begin{minipage}{0.3\textwidth}
\centering
\includegraphics[width=\linewidth]{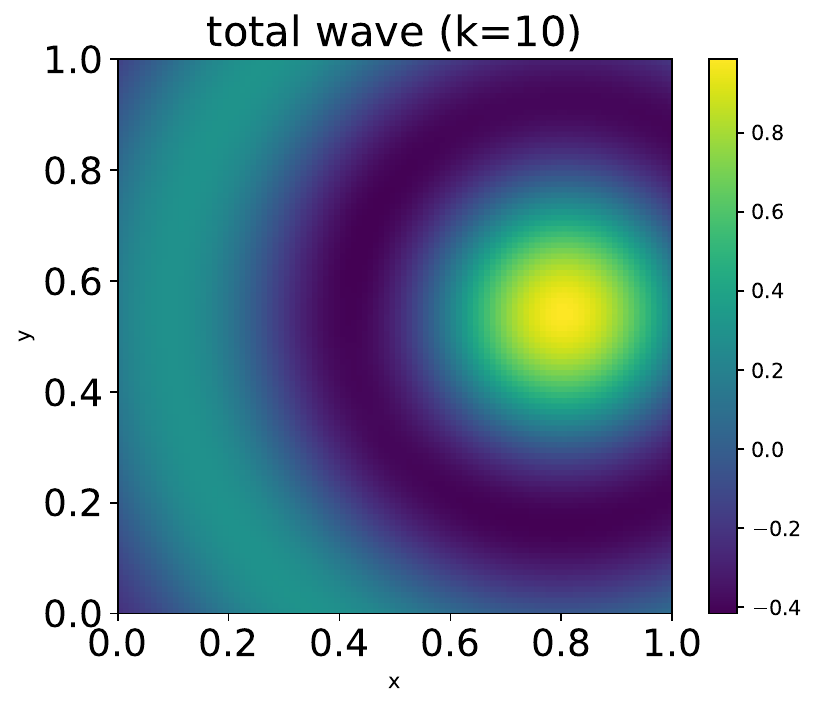}
\\ (a) $k=10$
\end{minipage}
\begin{minipage}{0.3\textwidth}
\centering
\includegraphics[width=\linewidth]{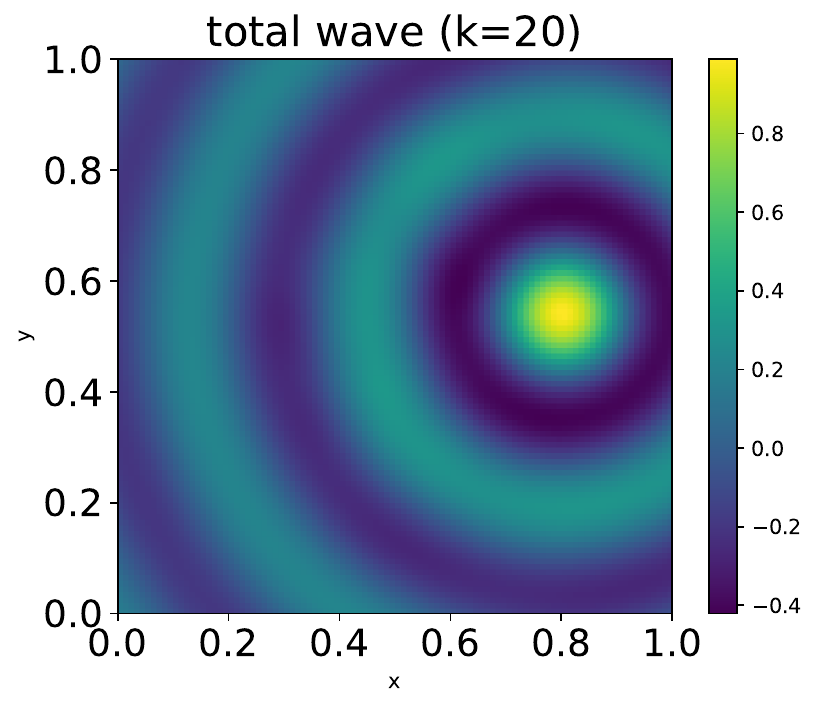}
\\ (b) $k=20$
\end{minipage}
\begin{minipage}{0.3\textwidth}
\centering
\includegraphics[width=\linewidth]{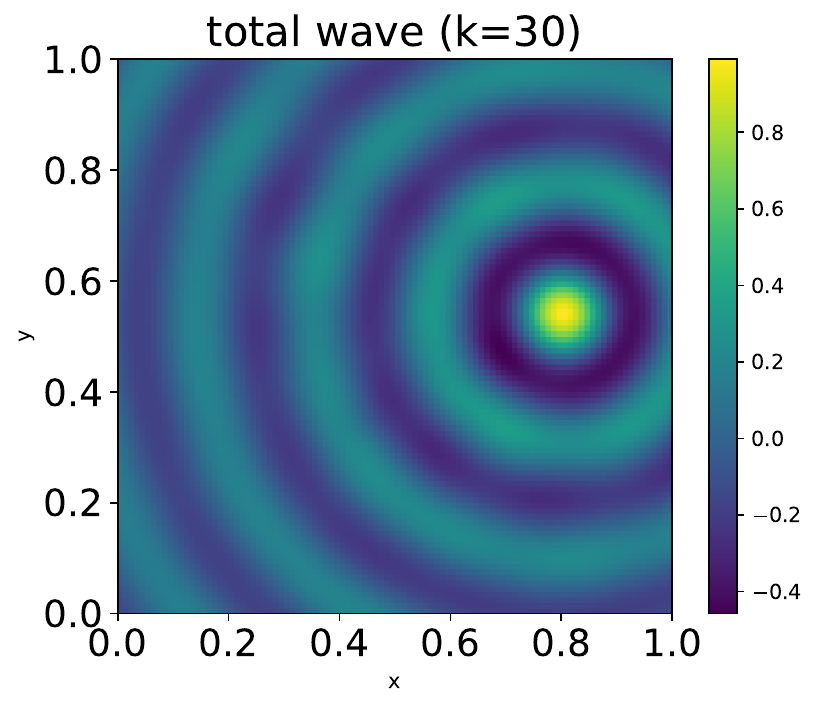}
\\ (c) $k=30$
\end{minipage}

\vspace{0.5em}

\begin{minipage}{0.3\textwidth}
\centering
\includegraphics[width=\linewidth]{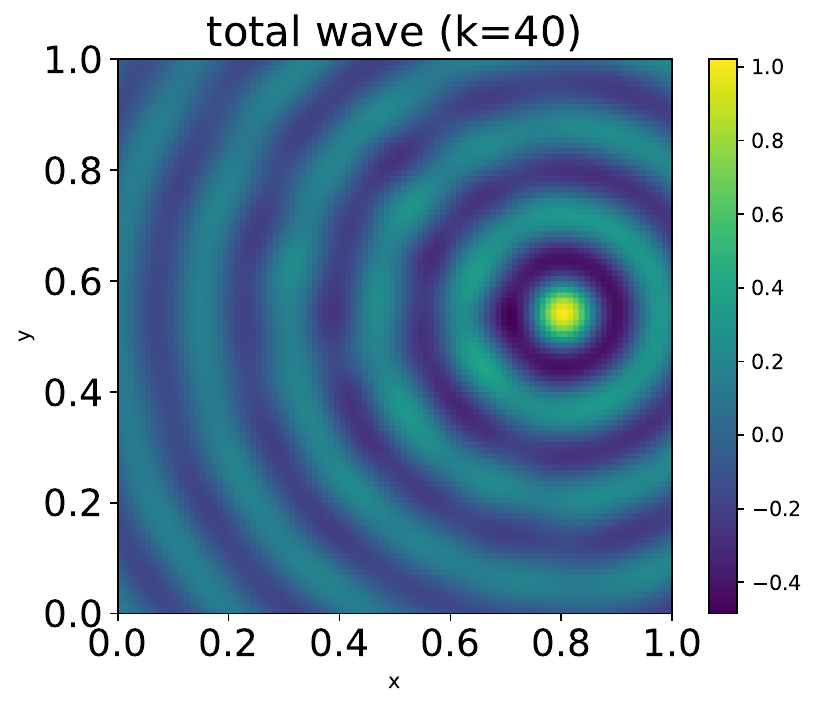}
\\ (d) $k=40$
\end{minipage}
\begin{minipage}{0.3\textwidth}
\centering
\includegraphics[width=\linewidth]{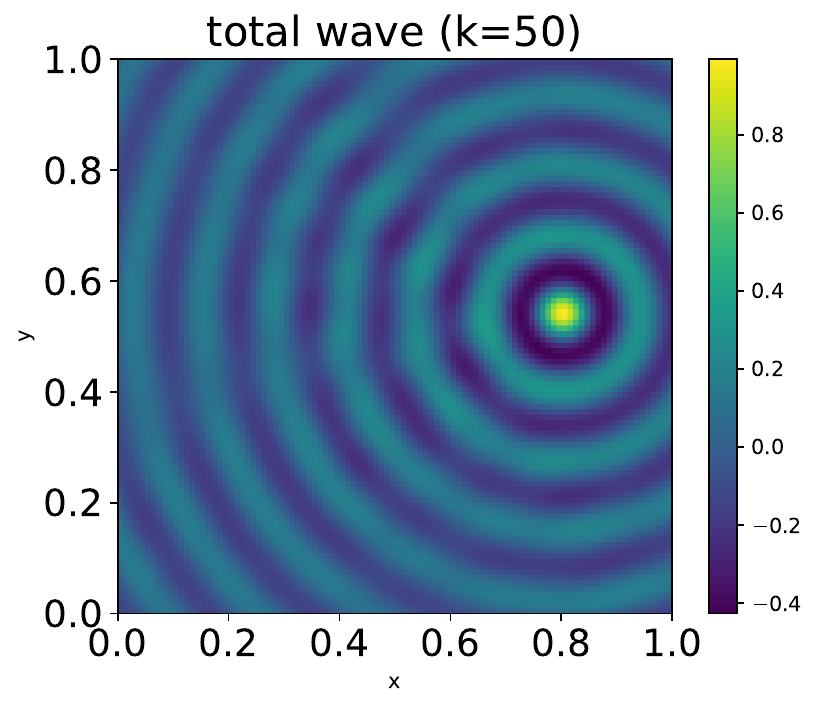}
\\ (e) $k=50$
\end{minipage}
\begin{minipage}{0.3\textwidth}
\centering
\includegraphics[width=\linewidth]{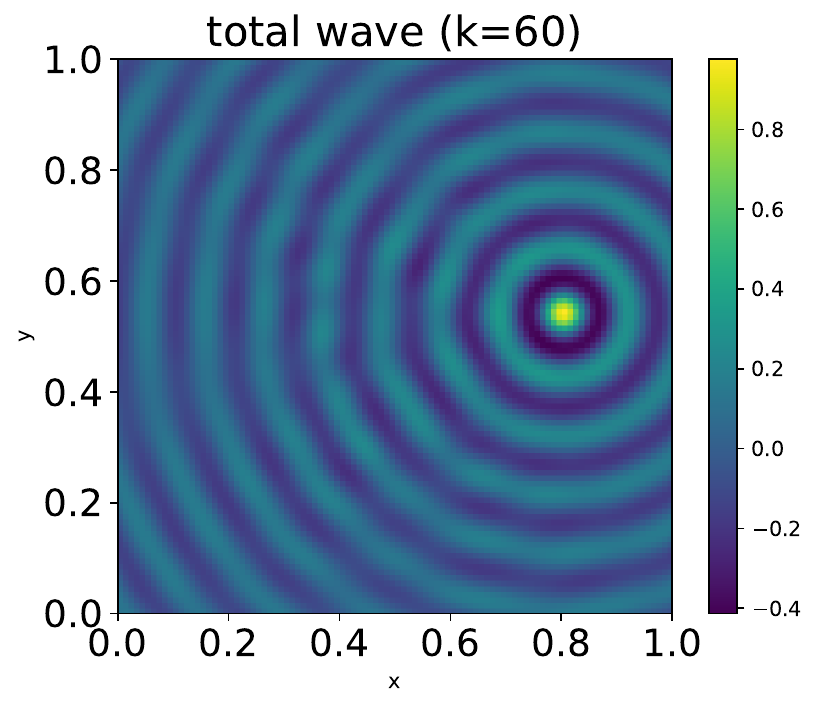}
\\ (f) $k=60$
\end{minipage}

\caption{Representative wavefields in the dataset for different wavenumbers.}
\label{fig:dataset_2d_wave}
\end{figure}

\noindent\textbf{Data generation.}
The source term is modeled as an ideal point source, $f(\bm x)=\delta(\bm x-\bm x_s).$
Synthetic datasets are generated by solving the Helmholtz equation~\eqref{helm_2d_eq} using a high-order finite element method, with relative numerical errors of approximately $10^{-3}$. 
Figure~\ref{fig:dataset_2d_wave} illustrates representative wavefields corresponding to different wavenumbers. The wavenumber $k$ is varied in the range $k \in [10,60]$ to assess performance across different frequency regimes. 
As the wavenumber $k$ increases, the wavefields exhibit increasingly complex oscillatory patterns and finer spatial structures, reflecting the increasing dominance of high-frequency oscillations.

\begin{figure}[t!]
\centering
\begin{minipage}{0.45\textwidth}
\centering
\includegraphics[width=\linewidth]{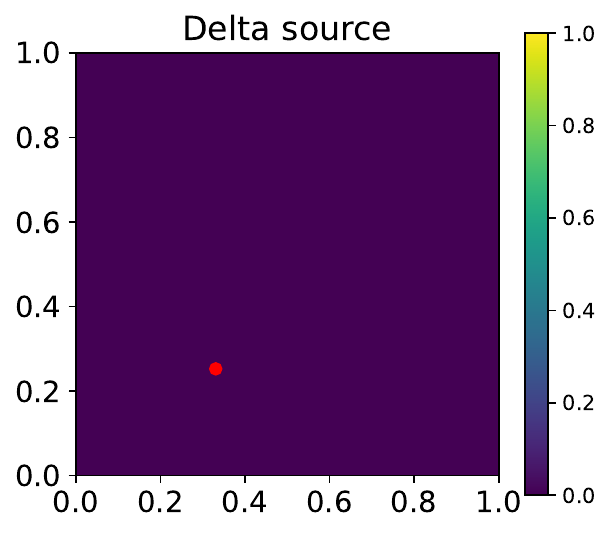}
\\ (a) Point source
\end{minipage}
\hspace{0.05\textwidth}
\begin{minipage}{0.45\textwidth}
\centering
\includegraphics[width=\linewidth]{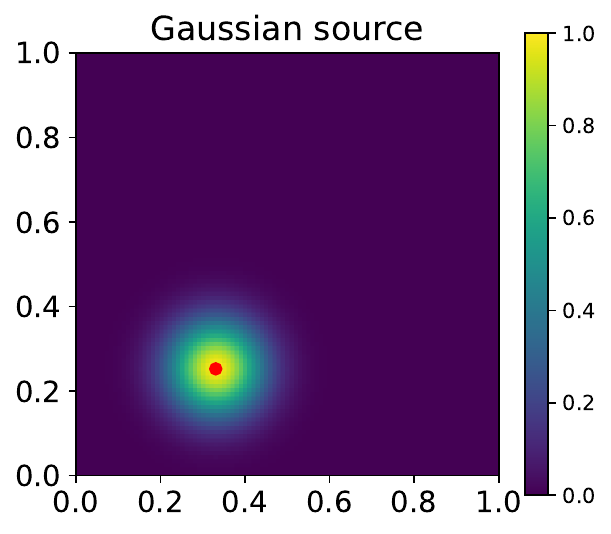}
\\ (b) Gaussian approximation
\end{minipage}
\caption{Comparison between point source and its Gaussian approximation in 2D.}
\label{fig:source_approx}
\end{figure}

\noindent\textbf{Source approximation.}
For neural operator learning, the delta source is approximated by a Gaussian function:
\begin{equation}
f(\bm x) \approx \exp\left(-\frac{\|\bm x - \bm x_s\|^2}{2\sigma_\epsilon^2}\right).
\end{equation}
This smoothing is applied only for the neural operator input, while the ground-truth data are generated using exact point sources.
A visual comparison between the ideal point source and its Gaussian approximation is shown in Figure~\ref{fig:source_approx}. The Gaussian representation preserves the localized nature of the source while providing a smooth and numerically stable input for the neural operator.

We evaluate both forward and inverse performance to assess the accuracy, stability, and high-frequency resolution capability of the proposed approach. In the forward setting, we measure the accuracy of predicted wavefields over a range of wavenumbers using FNO and MscaleFNO. In the inverse setting, we reconstruct the medium perturbation \(\alpha(\bm x)\) from observed data using the proposed plug-and-play diffusion-prior-guided inversion framework.

\subsection{Validation of the learned forward operator}
\label{sec:forward}
In this section,
we learn a family of forward operators $\mathcal{K}_k$ parameterized by the wavenumber $k$ using FNO/MscaleFNO, which map the medium and source to the corresponding wavefield:
\begin{equation}
\mathcal{K}_k(\alpha(x,y), f(x,y), x, y) \;\rightarrow\; u(x,y).
\end{equation}
For each wavenumber $k$, an independent dataset of $M \times R = 3000$ samples is constructed, consisting of triplets

\begin{equation}
\{\,\alpha_r^k(\bm x),\,f_m^k(\bm x),\,u_{m}^{r,k}(\bm x)\,\}_{m,r=1}^{M,R},
\end{equation}
where $\alpha_r^k$ denotes a medium realization, $f_m^k$ a source term, and $u_{m}^{r,k}$ the corresponding wavefield solution. Each dataset is split independently into training, validation, and test sets with a ratio of $8:1:1$.

To ensure fair comparison, both models have approximately the same number of parameters. The models utilize 20 Fourier modes along each spatial dimension and are trained using Adam optimizer \cite{kingma2015adam} with initial learning rate $\eta_{\max}=0.001$, and cosine annealing learning rate schedule \cite{loshchilov2017sgdr}. The normal FNO architecture equipped with 32 channels and 4 Fourier layers, comprising  \textbf{3,294,114} parameters. The parallel MscaleFNO framework with $4$ sub-networks, each branch operating with 16 channels and 4 Fourier layers, comprising \textbf{3,295,056} parameters. More implementation details can be found in \cite{You2026}.

\begin{figure}[htbp]
\centering
\includegraphics[width=0.6\textwidth]{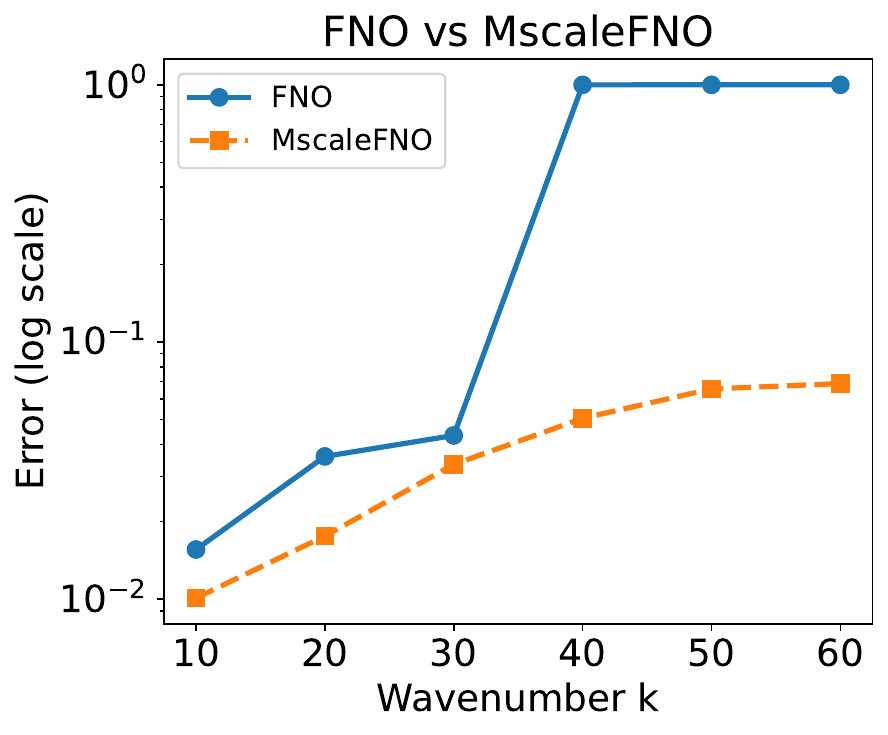}
\caption{Forward operator relative error versus wavenumber $k$ for FNO and MscaleFNO.}
\label{fig:error_vs_k_2d}
\end{figure}

A direct comparison is provided in Figure~\ref{fig:error_vs_k_2d}, which reports the relative prediction error as a function of the wavenumber $k$. Let $\bm{c}=\{c_i\}_{i=1}^4$ and $\bm{\gamma}=\{\gamma_i\}_{i=1}^4$ denote the 
scaling factors and weights of MscaleFNO, respectively. To better capture oscillatory 
solutions at different frequencies, MscaleFNO employs adaptive scale settings $\bm{c}$ 
for varying wavenumbers $k$: for $k=10,\;20$, $\bm{c} = \{10,\;50,\;100,\;200\}$; 
for $k=30,\;40$, $\bm{c} = \{10,\;50,\;200,\;1000\}$; and for $k=50,\;60$, 
$\bm{c} = \{10,\;100,\;500,\;2000\}$, with initial weights $\gamma_i=1,\;i=1,\dots,4$. The error of FNO increases significantly with $k$, while MscaleFNO remains stable and accurate, demonstrating strong robustness in the high-frequency regime.

These observations are consistent with the known spectral bias of neural operators, where low-frequency components are learned more efficiently than high-frequency ones. As a result, standard FNO struggles to approximate highly oscillatory solutions. The multiscale design of MscaleFNO mitigates this limitation by enabling effective representation across multiple frequency bands, leading to improved accuracy for large $k$. The improved forward accuracy is particularly important for inverse problems, where reconstruction quality critically depends on the fidelity of the forward operator. The results indicate that MscaleFNO provides a reliable surrogate model for high-frequency wave propagation, thereby serving as an accurate and efficient emulator for the subsequent inversion experiments.

\subsection{Inverse reconstruction}

We present inversion results for recovering the medium perturbation from noisy wavefield observations.
The inversion follows the diffusion-regularized framework in Algorithm~\ref{alg:pnp_diffusion}, where
the trained MscaleFNO surrogate is used as the differentiable forward operator and combined with an
EDM-based diffusion prior to improve reconstruction stability. The surrogate approximates the parameter-to-observation map from the medium perturbation and
source configuration to the wavefield measurements. This avoids repeatedly solving the Helmholtz
equation during iterative inversion and therefore substantially reduces the computational cost.

For each inversion experiment, we use $M=50$ source terms, with source locations distributed on a circle as described in Section~\ref{sec:experiment setup} and illustrated in Figure~\ref{fig:dataset_2d_medium}. The clean observation data are generated by solving the Helmholtz equation with a high-order finite element method, rather than by the learned surrogate model, to reduce the risk of inverse crime. The data are collected at $N=240$ receiver locations placed at grid points closest to the same observation circle used for the source distribution.
The clean observations are contaminated by additive Gaussian noise 
applied independently to the real and imaginary parts:
$u^{\delta} = u + \eta^{\delta}$,
where $\operatorname{Re}\eta^{\delta}, \operatorname{Im}\eta^{\delta} 
\sim \mathcal{N}(0,\sigma^2)$ with 
$\sigma = \sqrt{2}\delta\cdot\max(|u|)$.
We set $\delta = 10\%$ throughout and write $Y^\delta$ for the 
resulting noisy data.

The main inversion parameters used in Algorithm~\ref{alg:pnp_diffusion} are specified as follows. The total number of optimization iterations is set to $I = 1000$, and the Adam learning rate is chosen as $\eta = 0.001$. The EDM prior is applied periodically every $T = 30$ iterations, but starting after the first 200 iterations. This delay gives the optimizer room to explore the data-fidelity landscape and establish a coarse approximation before the denoiser begins pulling the estimates toward the learned manifold. Each EDM sampling uses $S = 4$ denoising steps.
The noise schedule used in the EDM sampler is updated during the inversion process. Since the EDM noise sequence $\{\sigma_s\}_{s=0}^{S}$ is arranged in descending order from a maximum value to a minimum value, we denote the largest and smallest noise levels at a given EDM update by $\sigma_{\max}^{\mathrm{now}}=\sigma_S$ and $\sigma_{\min}=\sigma_0$, respectively. In all experiments, the minimum noise level is fixed as
$\sigma_{\min} = 0.002$.
The maximum noise level is gradually decreased as the inversion iteration proceeds. Let $i$ denote the current inversion iteration and $I$ be the total number of iterations. We define the normalized iteration variable by $\tau_i = (i+1)/I$. 
Then the current maximum noise level is linearly interpolated between the initial value $\sigma_{\max}^{\mathrm{init}} = 0.3$ and the final value $\sigma_{\max}^{\mathrm{final}} = 0.01$:
$
\sigma_{\max}^{\mathrm{now}} 
=
(1-\tau_i)\sigma_{\max}^{\mathrm{init}}
+
\tau_i \sigma_{\max}^{\mathrm{final}}.
$
Therefore, as the inversion progresses, the EDM sampler gradually shifts from stronger prior correction to weaker fine-scale refinement.

For each EDM application, the sigma sequence is generated using the Karras schedule. Given the EDM sampling step number $S$, the $s$th sigma value, where $s=0,1,\ldots,S$, is computed as
$t_s = \frac{s}{S}$, and
\[
\sigma_s
=
\left(
\left(\sigma_{\max}^{\mathrm{now}}\right)^{1/\rho}
+
t_s
\left(
\sigma_{\min}^{1/\rho}
-
\left(\sigma_{\max}^{\mathrm{now}}\right)^{1/\rho}
\right)
\right)^\rho .
\]
Here we set $\rho = 7$.
This schedule performs a linear interpolation in the transformed $\sigma^{1/\rho}$ space and then maps the values back to the original sigma space. As a result, the sequence decreases from $\sigma_{\max}^{\mathrm{now}}$ to $\sigma_{\min}$, while allocating more sampling resolution to the low-noise region, which is important for fine-scale reconstruction.

\begin{figure}[t!]
\centering
\includegraphics[width=0.6\textwidth]{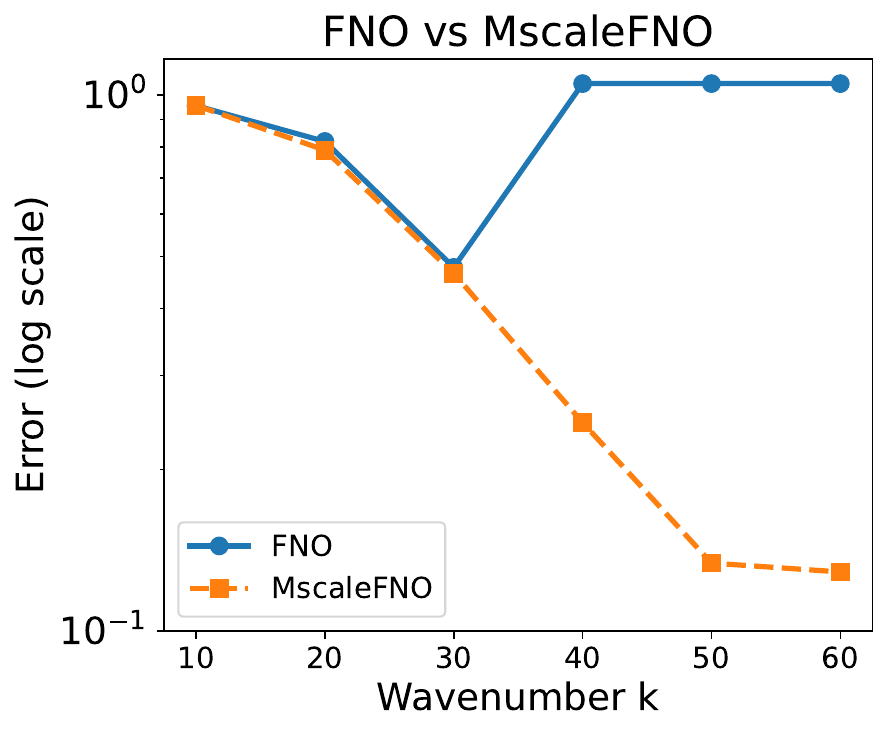}
\caption{Inversion error across different wavenumbers $k$ under $10\%$ noise.}
\label{fig:inv2d_error_compare}
\end{figure}

Figure~\ref{fig:inv2d_error_compare} compares the reconstruction errors of FNO and MscaleFNO across different wavenumbers. As the wavenumber $k$ increases, the performance of FNO deteriorates significantly, with consistently high errors in the high-frequency regime. In contrast, MscaleFNO effectively utilizes high-frequency information and achieves substantially lower reconstruction errors across all tested cases. This is consistent with the forward operator validation in Section~\ref{sec:forward}, where FNO also struggled at high frequencies.

\begin{figure}[htbp]
\centering

\begin{subfigure}{0.9\textwidth}
    \centering
    \includegraphics[width=0.44\linewidth,height=0.13\textheight]{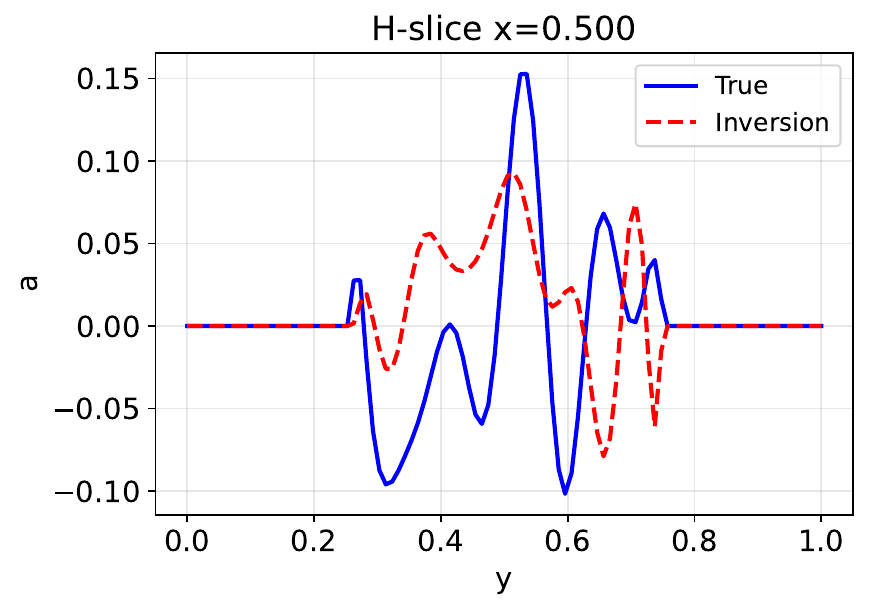}
    \includegraphics[width=0.44\linewidth,height=0.13\textheight]{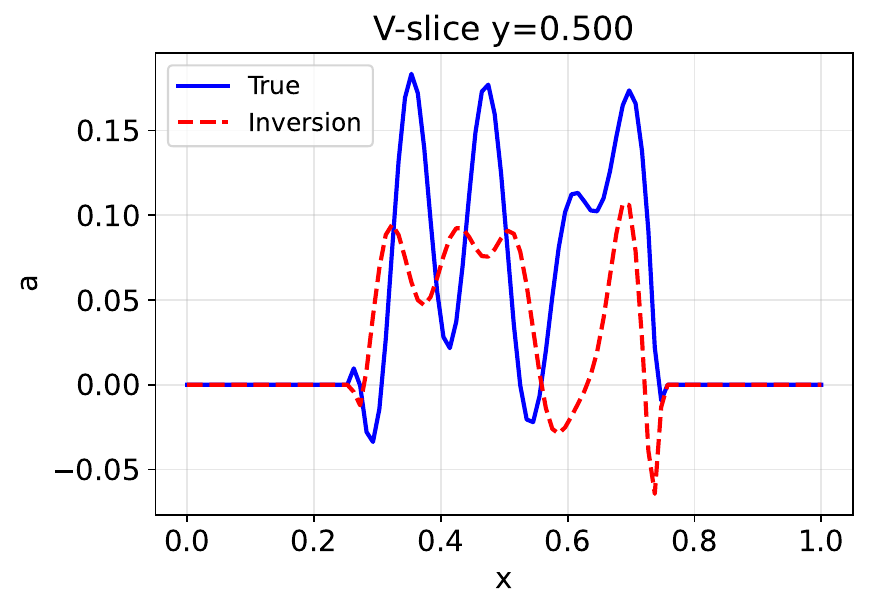}
\caption{$k=10,\ \mathrm{err}=0.956$}
\end{subfigure}

\vspace{0.1em}

\begin{subfigure}{0.9\textwidth}
    \centering
    \includegraphics[width=0.44\linewidth,height=0.13\textheight]{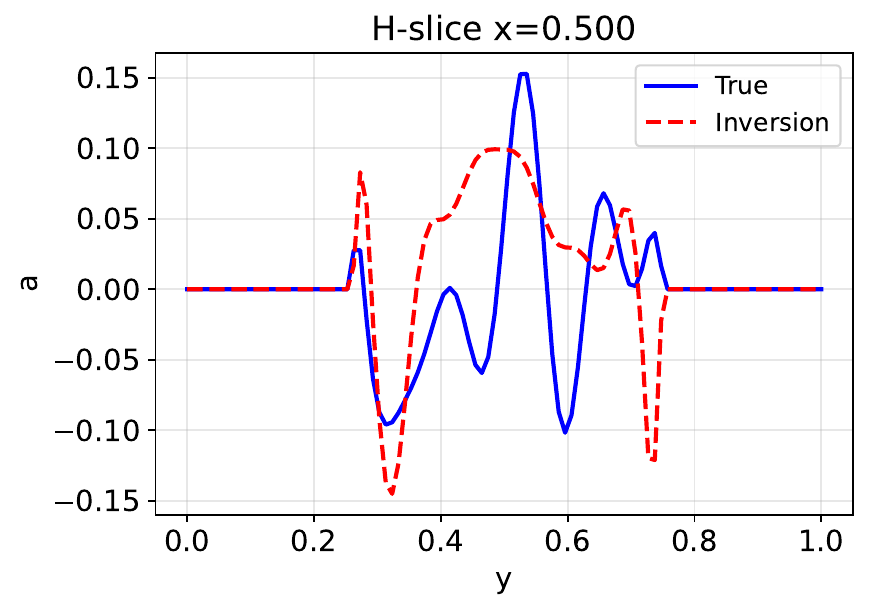}
    \includegraphics[width=0.44\linewidth,height=0.13\textheight]{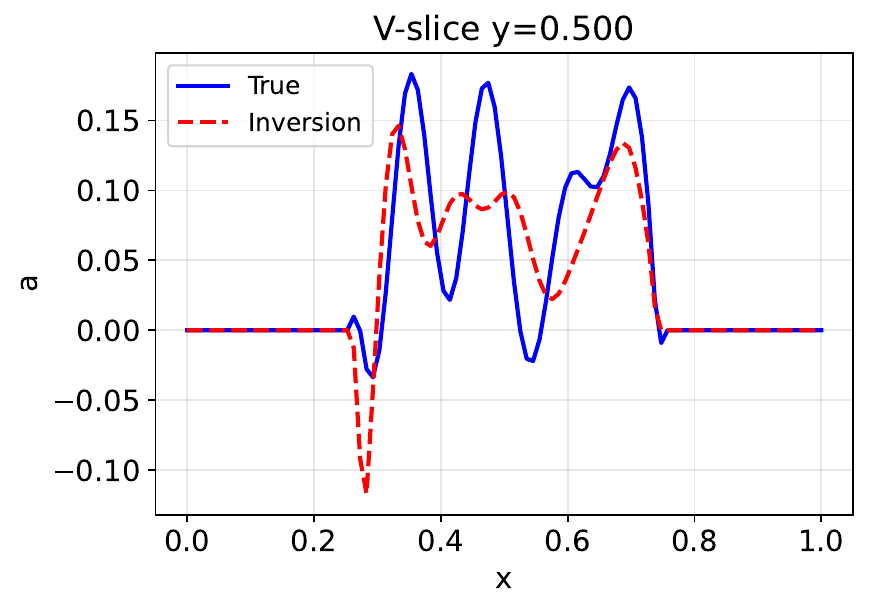}
\caption{$k=20,\ \mathrm{err}=0.784$}
\end{subfigure}

\vspace{0.1em}

\begin{subfigure}{0.9\textwidth}
    \centering
    \includegraphics[width=0.44\linewidth,height=0.13\textheight]{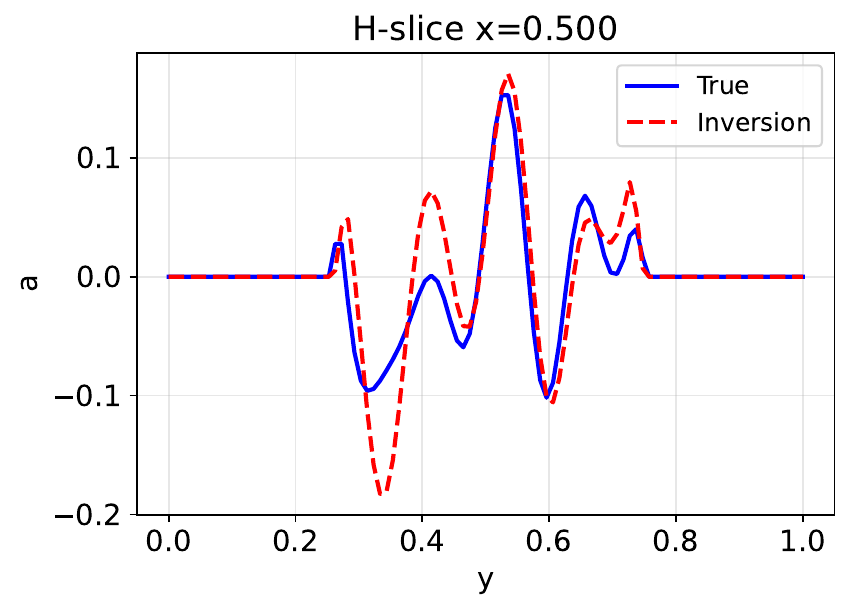}
    \includegraphics[width=0.44\linewidth,height=0.13\textheight]{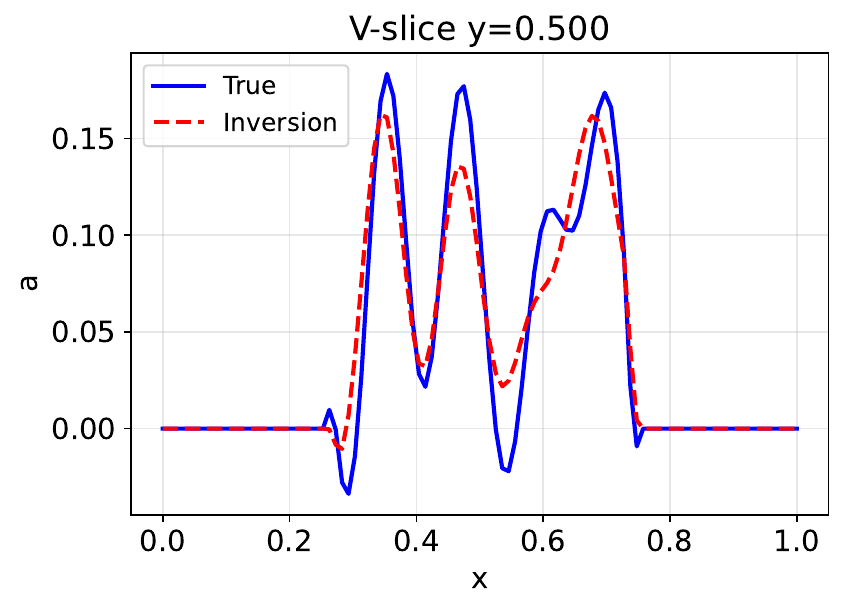}
\caption{$k=30,\ \mathrm{err}=0.465$}
\end{subfigure}

\vspace{0.1em}

\begin{subfigure}{0.9\textwidth}
    \centering
    \includegraphics[width=0.44\linewidth,height=0.13\textheight]{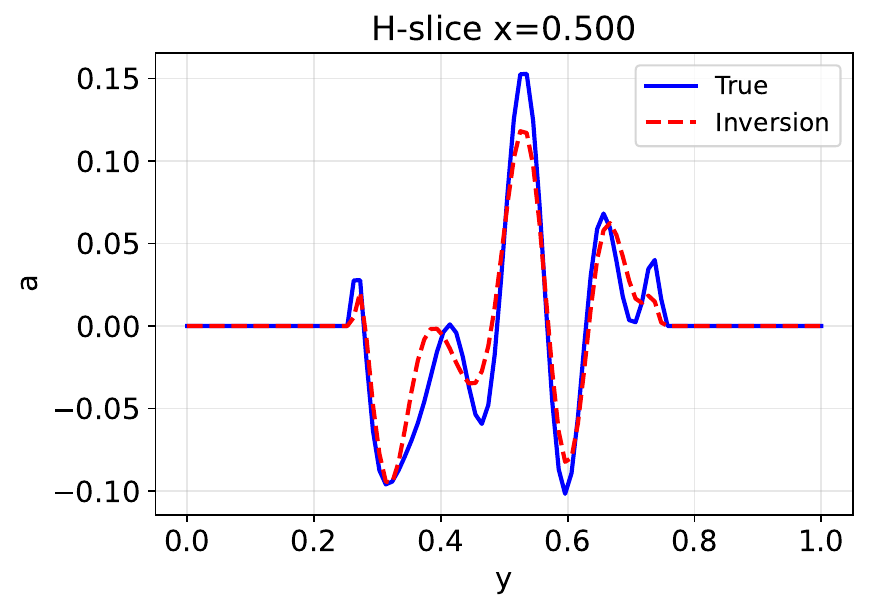}
    \includegraphics[width=0.44\linewidth,height=0.13\textheight]{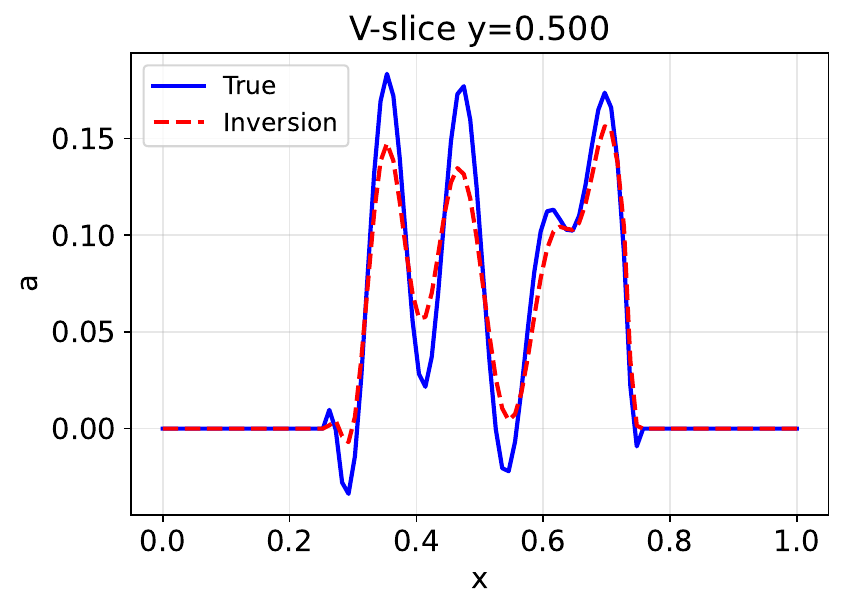}
\caption{$k=40,\ \mathrm{err}=0.245$}
\end{subfigure}

\vspace{0.1em}

\begin{subfigure}{0.9\textwidth}
    \centering
    \includegraphics[width=0.44\linewidth,height=0.13\textheight]{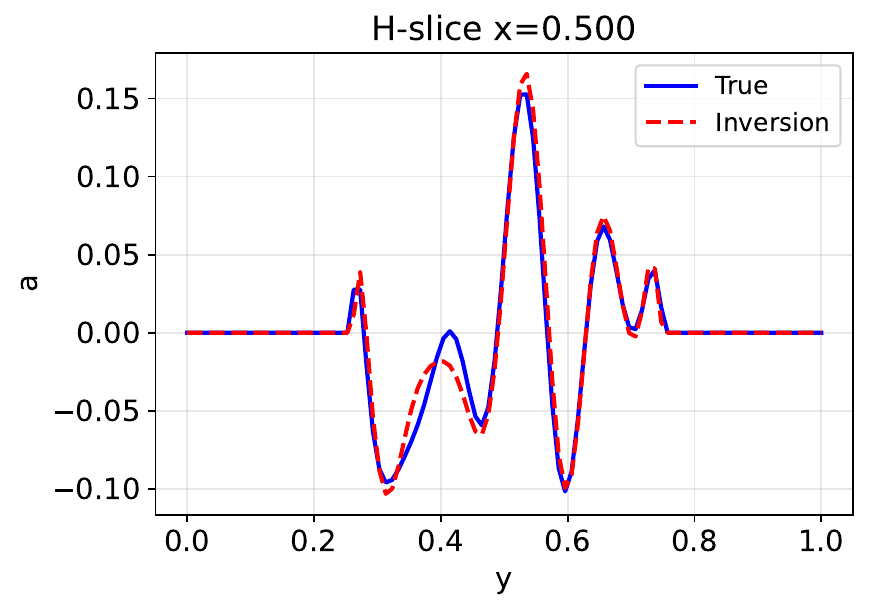}
    \includegraphics[width=0.44\linewidth,height=0.13\textheight]{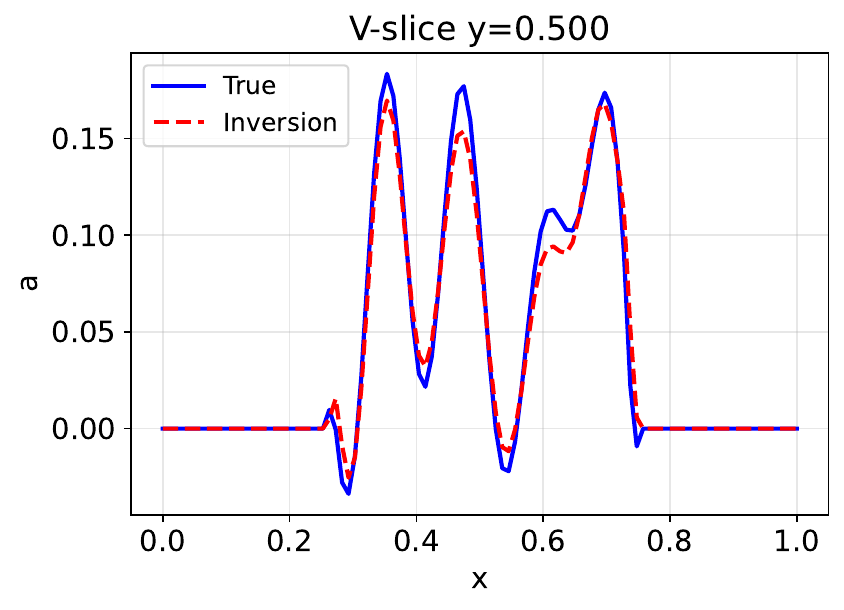}
\caption{$k=50,\ \mathrm{err}=0.134$}
\end{subfigure}

\vspace{0.1em}

\begin{subfigure}{0.9\textwidth}
    \centering
    \includegraphics[width=0.44\linewidth,height=0.13\textheight]{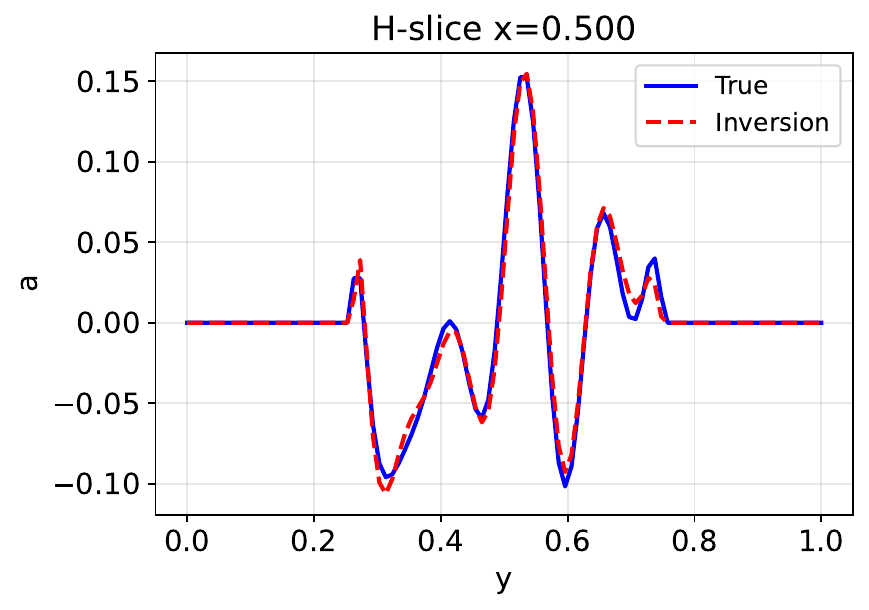}
    \includegraphics[width=0.44\linewidth,height=0.13\textheight]{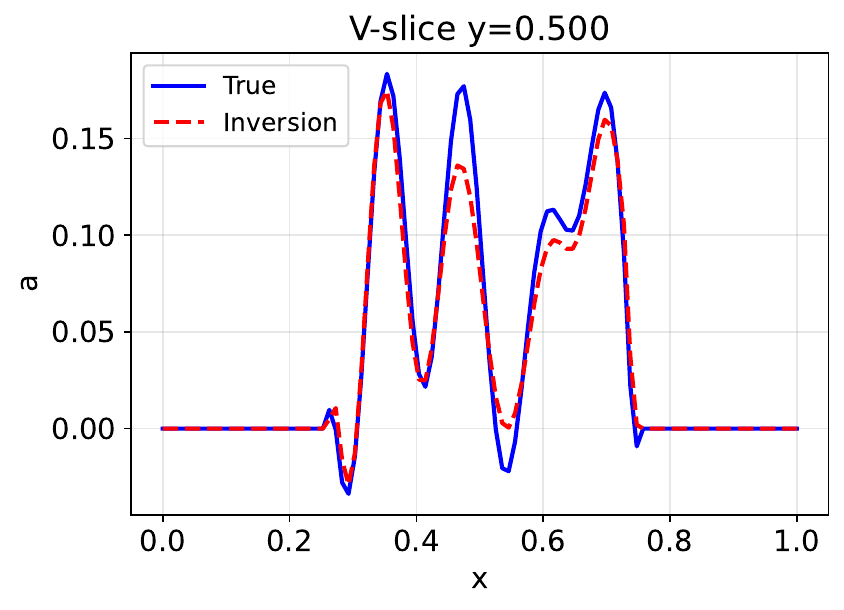}
\caption{$k=60,\ \mathrm{err}=0.129$}
\end{subfigure}

\caption{Cross-sectional comparisons of reconstruction results by MscaleFNO for different wavenumbers 
$k$ under $10\%$ noise. }
\label{fig:inv2d_reconstruction}
\end{figure}

Figure~\ref{fig:inv2d_reconstruction} provides cross-sectional comparisons of the reconstructed media for different values of $k$. As the wavenumber increases, the reconstruction quality improves markedly, with finer-scale structures becoming progressively better resolved. This behavior is consistent with the fact that higher-frequency waves carry richer spatial information, which enhances the identifiability of fine features in the inverse problem. Notably, the reconstruction error decreases from 0.956 at $k=10$ to 0.129 at $k=60$, confirming the benefit of high-frequency measurements.

\begin{figure}[htbp]
\centering
\includegraphics[width=0.55\textwidth]{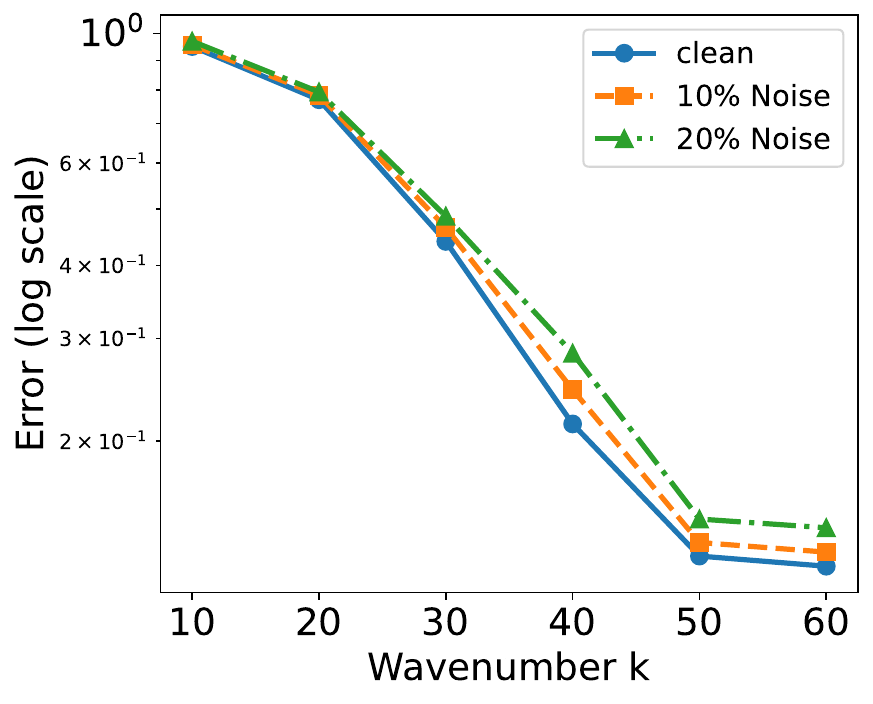}
\caption{Inversion error versus wavenumber $k$ under different noise levels.}
\label{fig:inv2d_noise}
\end{figure}

\noindent\textbf{Inversion robustness under noisy observations}
Figure~\ref{fig:inv2d_noise} shows the reconstruction error under different noise levels. While the errors increase as the noise level grows, MscaleFNO remains stable across all tested wavenumbers and maintains consistent performance even in the presence of significant noise. This demonstrates the robustness of the proposed approach for practical inverse problems with noisy measurements.

\noindent\textbf{Limited-aperture configuration} In practical settings, data acquisition is often confined to a limited angular range by physical and experimental constraints. This incomplete coverage, where sources and receivers span only a subset of directions, introduces additional ill-posedness.
We consider a two-dimensional setting where sources and receivers are distributed along a circular observation curve surrounding the medium, but confined to a finite angular interval. The objective is to systematically investigate how limited-aperture configurations affect the recovery of the medium $\alpha(\bm x)$. All reconstructions are performed using MscaleFNO with wavenumber $k=60$, and $10\%$ observational noise is added. We consider source apertures of $45^\circ$, $90^\circ$, $180^\circ$, and $360^\circ$, and independently vary the receiver apertures over the same range. We also test a ``diagonal'' receiver layout, where receivers are split into two equal sectors placed on opposite sides of the domain (e.g., diagonal $r=45^\circ$ means two $45^\circ$ sectors, totaling $90^\circ$).

\begin{figure}[t!]
\centering
\setlength{\tabcolsep}{2pt}

\begin{tabular}{cccc}

\begin{tabular}{c}
{\small $r=45^\circ, s=45^\circ$} \\
\includegraphics[height=2.8cm]{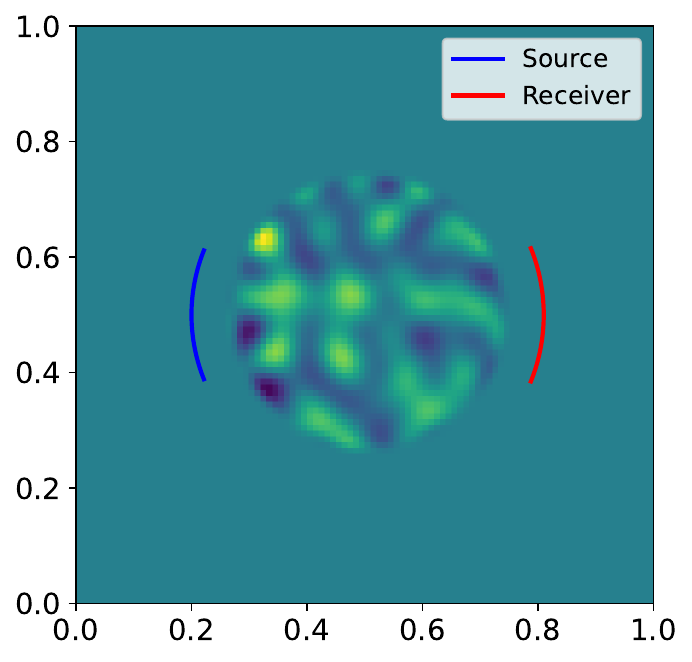} 
\end{tabular} &
\begin{tabular}{c}
{\small $r=45^\circ, s=90^\circ$} \\
\includegraphics[height=2.8cm]{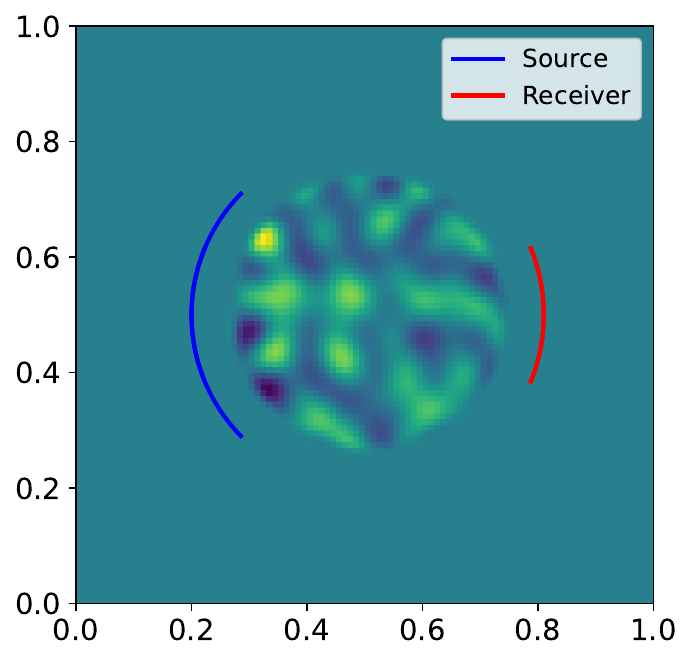} 
\end{tabular} &
\begin{tabular}{c}
{\small $r=45^\circ, s=180^\circ$} \\
\includegraphics[height=2.8cm]{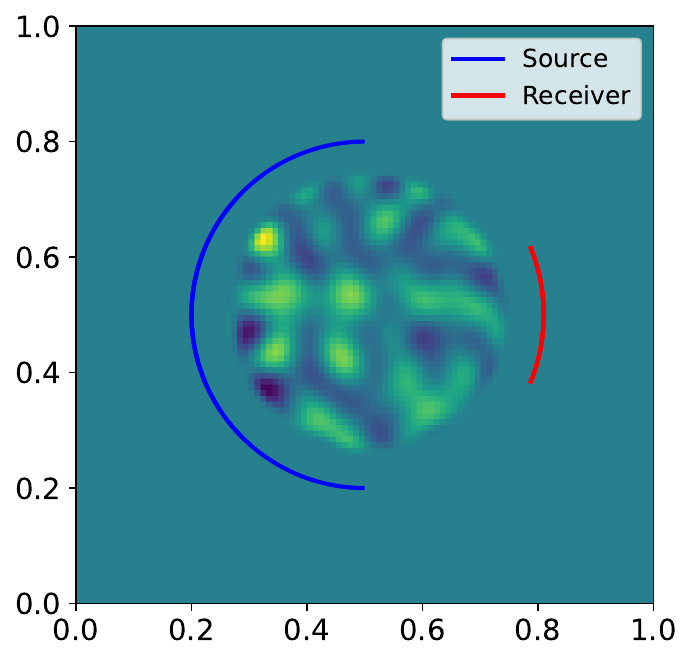} 
\end{tabular} &
\begin{tabular}{c}
{\small $r=45^\circ, s=360^\circ$} \\
\includegraphics[height=2.8cm]{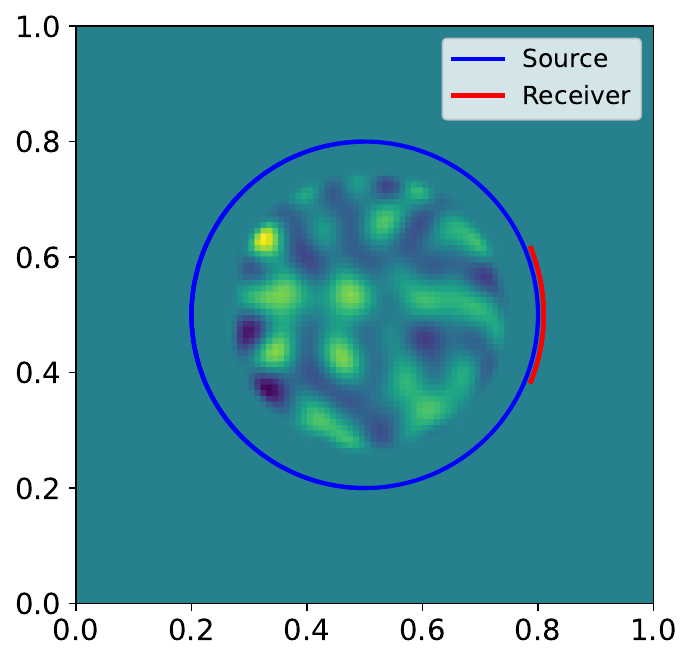} 
\end{tabular}
\\

\begin{tabular}{c}
{\small $r=90^\circ, s=45^\circ$} \\
\includegraphics[height=2.8cm]{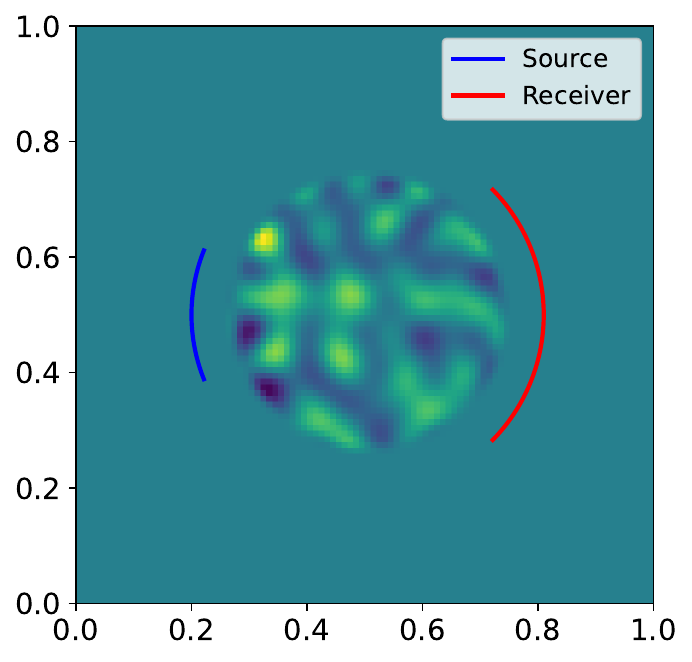} 
\end{tabular} &
\begin{tabular}{c}
{\small $r=90^\circ, s=90^\circ$} \\
\includegraphics[height=2.8cm]{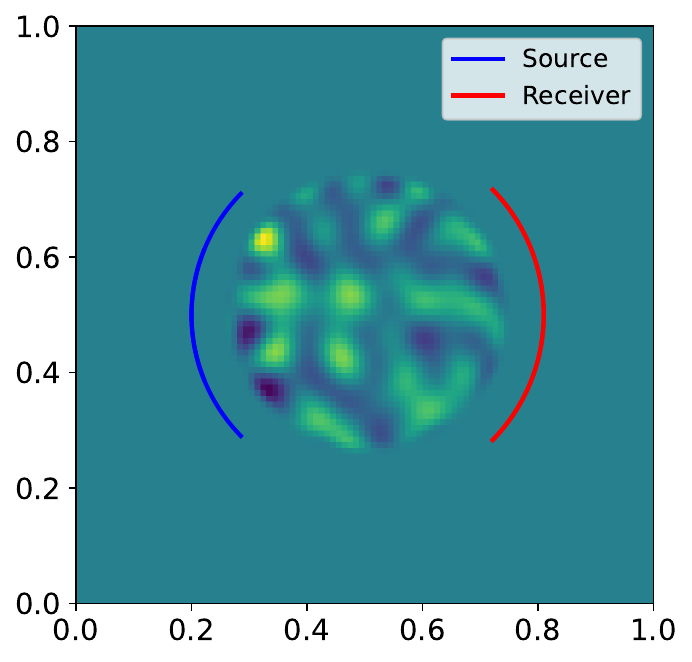} 
\end{tabular} &
\begin{tabular}{c}
{\small $r=90^\circ, s=180^\circ$} \\
\includegraphics[height=2.8cm]{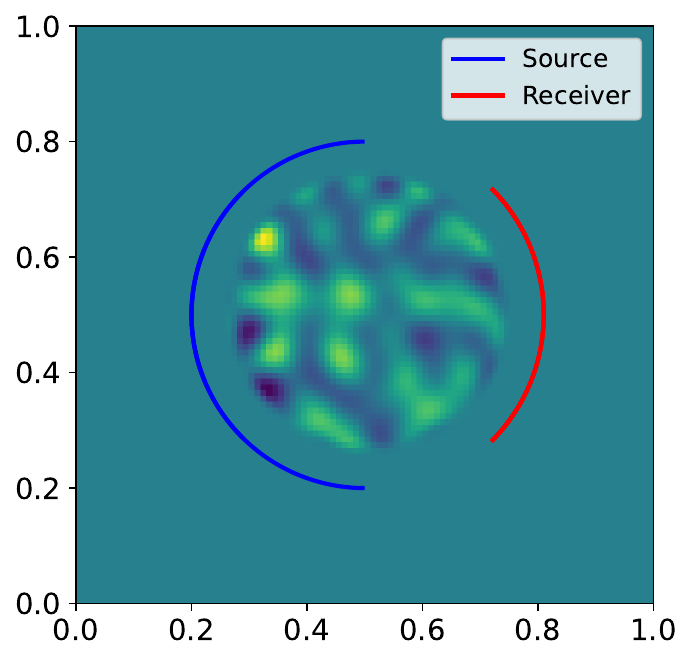} 
\end{tabular} &
\begin{tabular}{c}
{\small $r=90^\circ, s=360^\circ$} \\
\includegraphics[height=2.8cm]{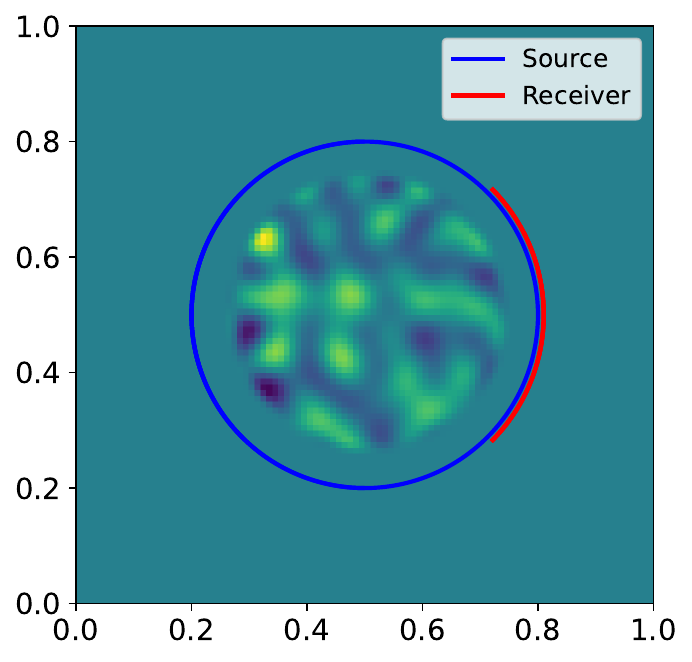} 
\end{tabular}
\\

\begin{tabular}{c}
{\small $r=180^\circ, s=45^\circ$} \\
\includegraphics[height=2.8cm]{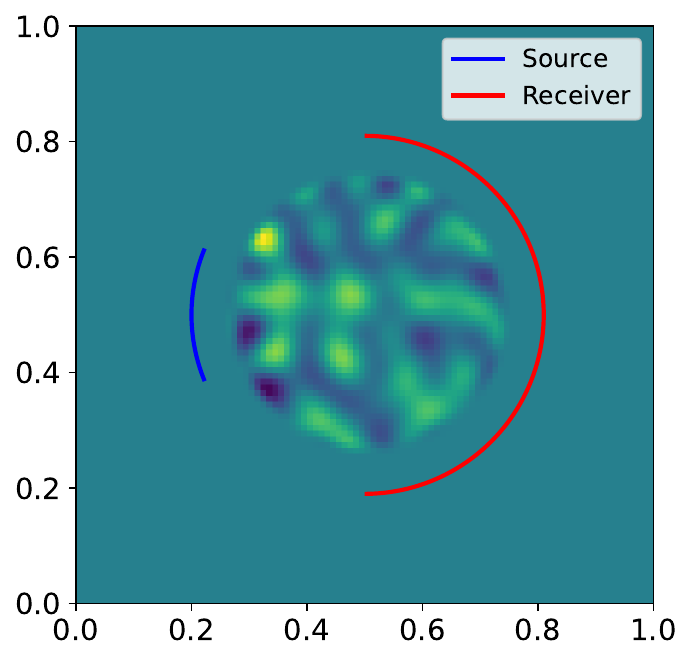} 
\end{tabular} &
\begin{tabular}{c}
{\small $r=180^\circ, s=90^\circ$} \\
\includegraphics[height=2.8cm]{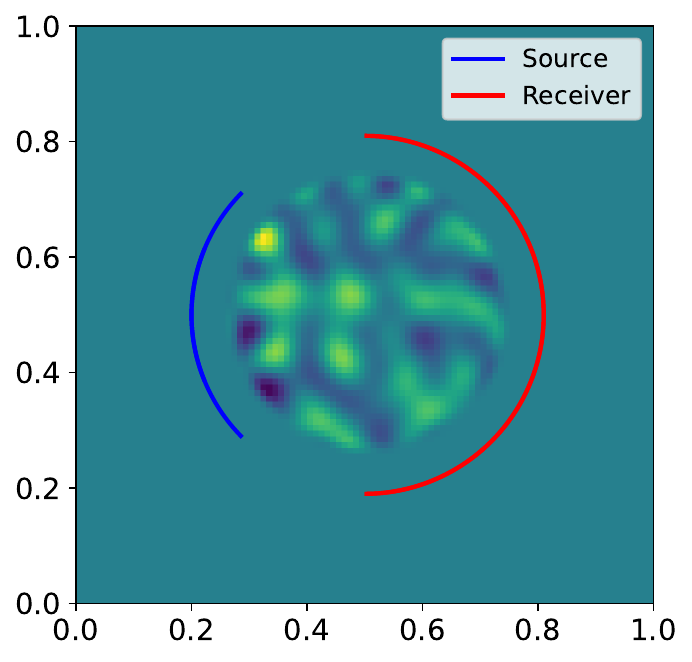} 
\end{tabular} &
\begin{tabular}{c}
{\small $r=180^\circ, s=180^\circ$} \\
\includegraphics[height=2.8cm]{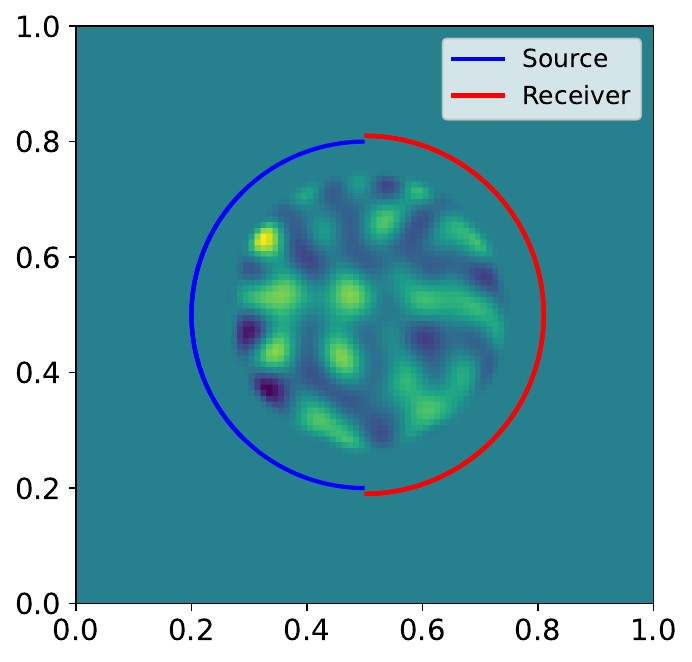}
\end{tabular} &
\begin{tabular}{c}
{\small $r=180^\circ, s=360^\circ$} \\
\includegraphics[height=2.8cm]{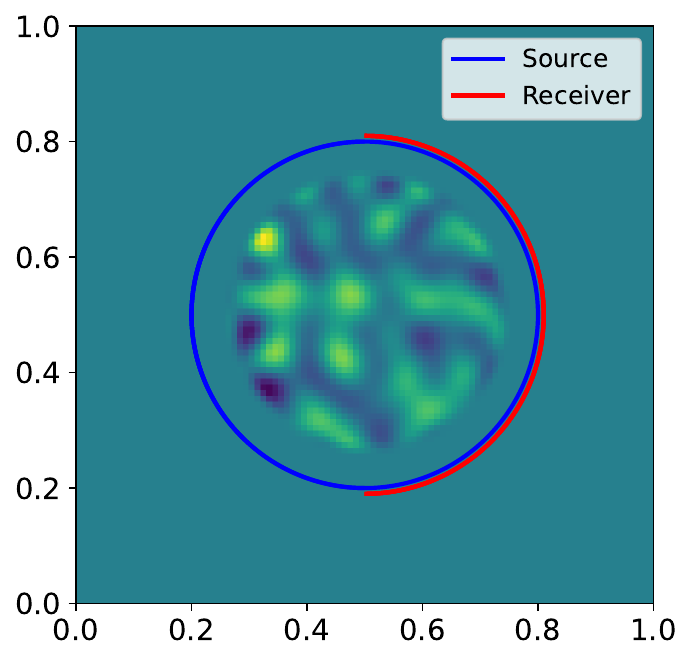}
\end{tabular}
\\

\begin{tabular}{c}
{\small $r=360^\circ, s=45^\circ$} \\
\includegraphics[height=2.8cm]{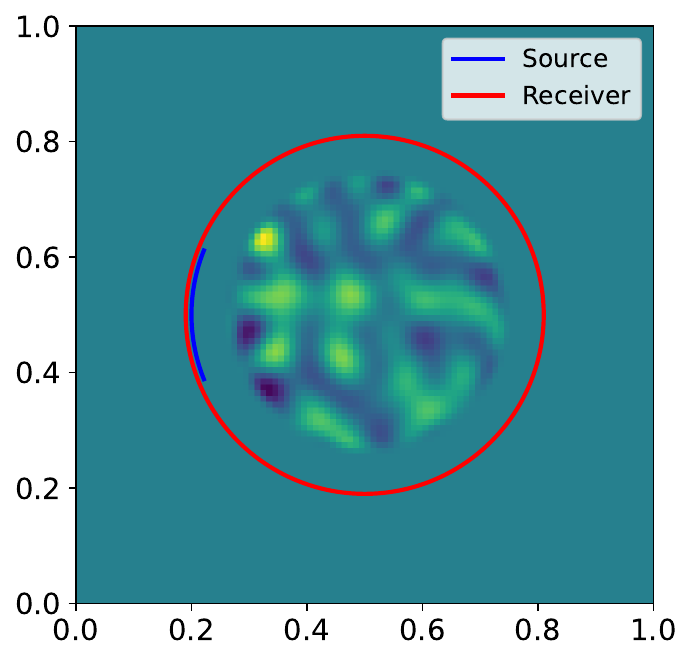} 
\end{tabular} &
\begin{tabular}{c}
{\small $r=360^\circ, s=90^\circ$} \\
\includegraphics[height=2.8cm]{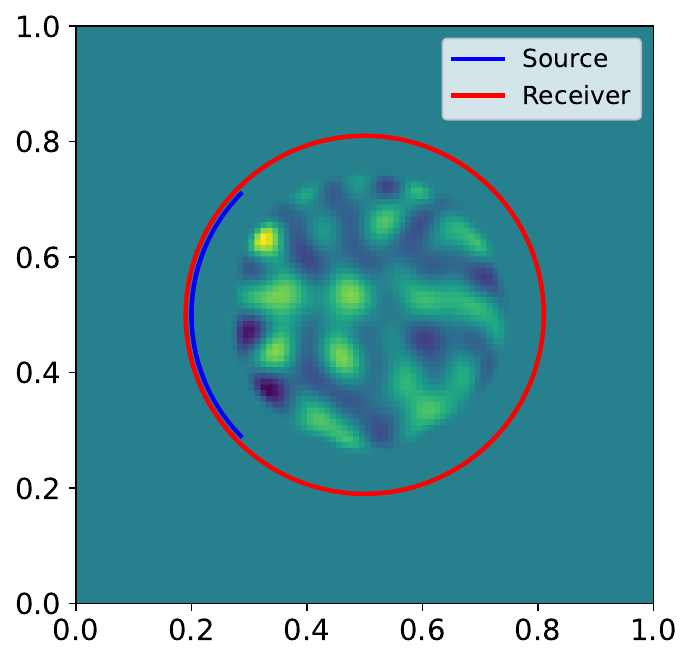} 
\end{tabular} &
\begin{tabular}{c}
{\small $r=360^\circ, s=180^\circ$} \\
\includegraphics[height=2.8cm]{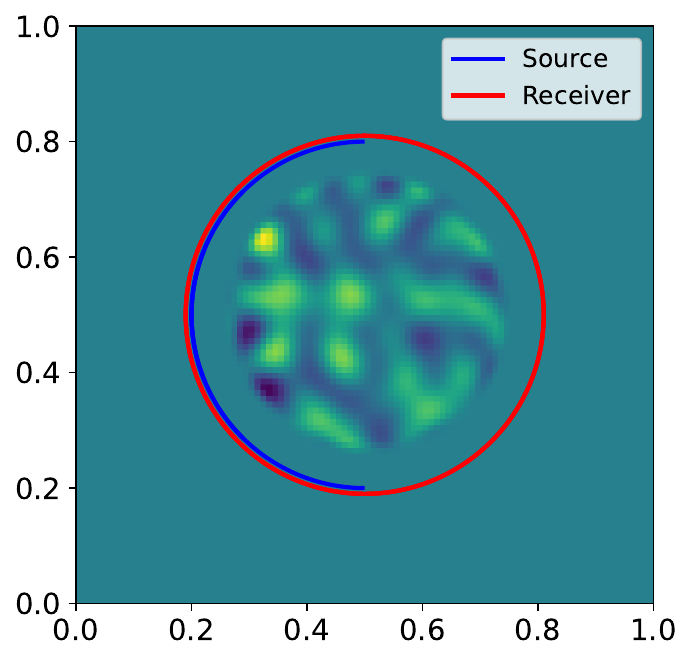} 
\end{tabular} &
\begin{tabular}{c}
{\small $r=360^\circ, s=360^\circ$} \\
\includegraphics[height=2.8cm]{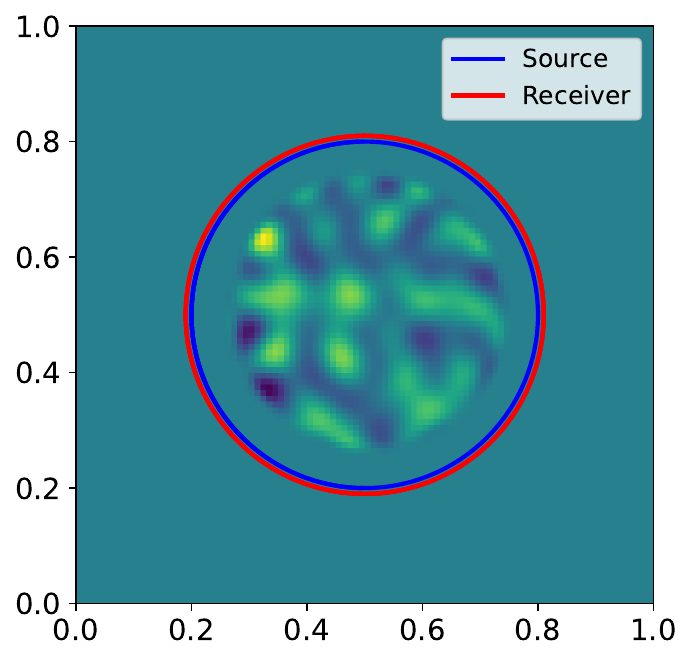} 
\end{tabular}

\end{tabular}

\centering
\begin{tabular}{cc}

\begin{tabular}{c}
{\small Diagonal $r=45^\circ, s=360^\circ$} \\
\includegraphics[height=2.8cm]{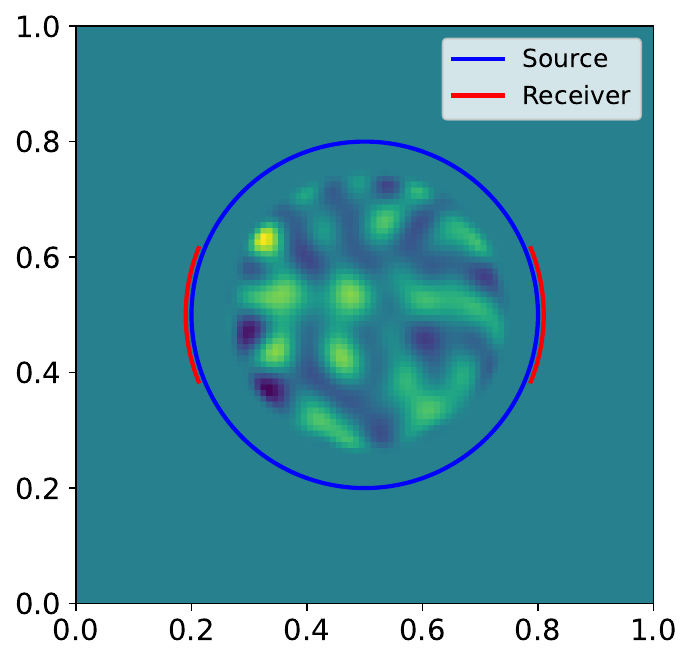} 
\end{tabular}

&

\begin{tabular}{c}
{\small Diagonal $r=90^\circ, s=360^\circ$} \\
\includegraphics[height=2.8cm]{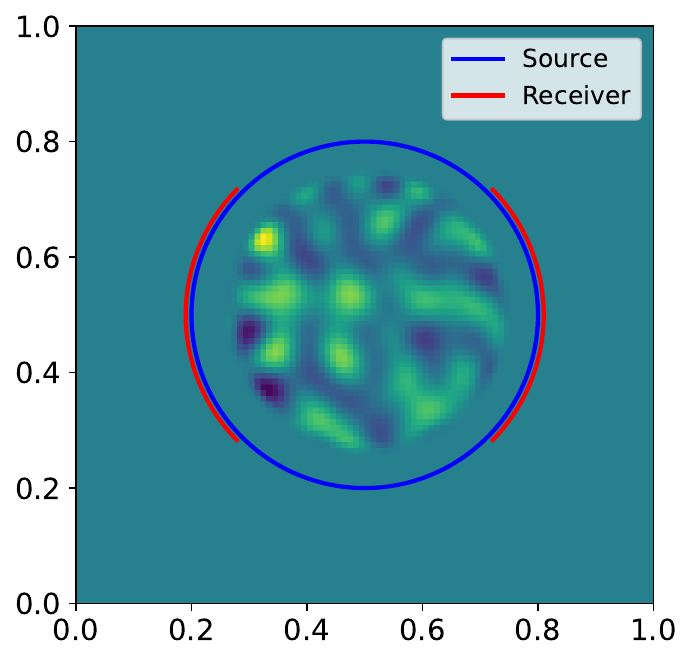} 
\end{tabular}

\end{tabular}
\caption{
Illustration of source (s) and receiver (r) configurations under different limited-aperture settings. Columns correspond to source apertures and rows correspond to receiver apertures, both ranging over $45^\circ$, $90^\circ$, $180^\circ$, and $360^\circ$. 
}
\label{fig:limited_aperture_configuration}
\end{figure}

Figure~\ref{fig:limited_aperture_configuration} illustrates the corresponding source and receiver configurations. Columns represent different source apertures and rows represent different receiver apertures. Each configuration corresponds to a specific combination of illumination and observation coverage.

\begin{table}[t!]
\centering
\caption{Relative reconstruction errors under different source ($s$) and receiver ($r$) aperture configurations.}
\label{tab:aperture_error}
\begin{tabular}{c|cccc}
\hline
$r \backslash s$ & $45^\circ$ & $90^\circ$ & $180^\circ$ & $360^\circ$ \\
\hline
$45^\circ$  & 1.000 & 0.947 & 0.744 & 0.476 \\
$90^\circ$  & 0.919 & 0.833 & 0.671 & 0.336 \\
$180^\circ$ & 0.833 & 0.684 & 0.435 & 0.143 \\
$360^\circ$ & 0.569 & 0.421 & 0.149 & \textbf{0.129} \\
\midrule
\multicolumn{5}{c}{\textit{Diagonal receiver configurations}} \\
\midrule
\multicolumn{3}{c}{$r=45^\circ,\, s=360^\circ$} & 
\multicolumn{2}{c}{$r=90^\circ,\, s=360^\circ$} \\
\multicolumn{3}{c}{0.319} & 
\multicolumn{2}{c}{0.157} \\
\bottomrule
\end{tabular}
\end{table}
Table~\ref{tab:aperture_error} reports the relative reconstruction errors across all source--receiver aperture combinations. The error decreases monotonically as either the source or receiver aperture increases, indicating improved reconstruction accuracy with enhanced angular coverage. The lowest errors are achieved in near full-aperture configurations (e.g., $360^\circ$), while severely limited apertures (e.g., $45^\circ$) result in significantly higher errors. 
Furthermore, the diagonal configurations highlight the benefit of spatial diversity. For instance, under full illumination ($s=360^\circ$), the diagonal $r=45^\circ$ configuration (which has a total receiver aperture of $90^\circ$) yields an error of 0.319, outperforming the contiguous $90^\circ$ receiver configuration (error = 0.336). This suggests that distributing the same angular coverage across opposite sides can further mitigate the ill-posedness compared to a single localized observation window.

\begin{figure}[t!]
\centering
\setlength{\tabcolsep}{1pt}

\begin{tabular}{cccc}

\begin{tabular}{c}
{\small $r=45^\circ, s=45^\circ$} \\
\includegraphics[height=2.4cm]{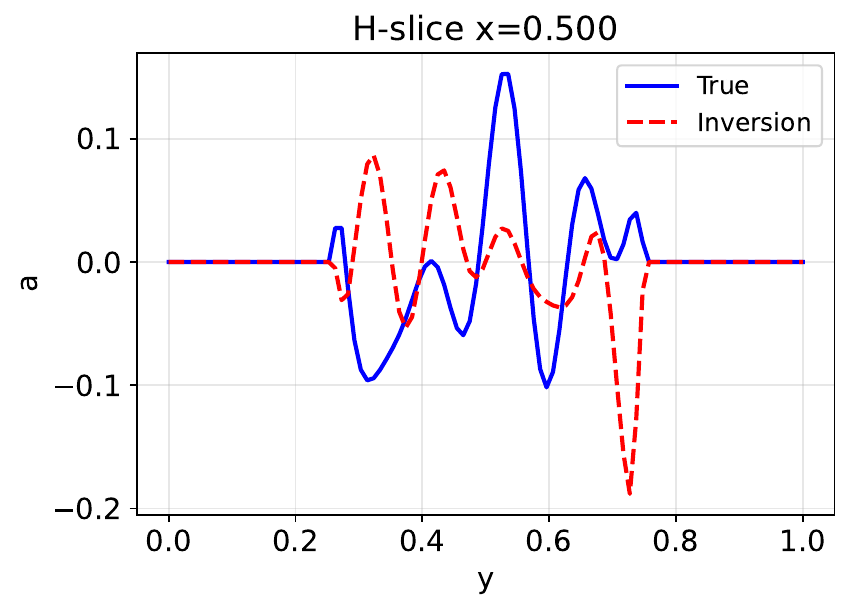} \\
{\scriptsize error=1.000}
\end{tabular} &
\begin{tabular}{c}
{\small $r=45^\circ, s=90^\circ$} \\
\includegraphics[height=2.4cm]{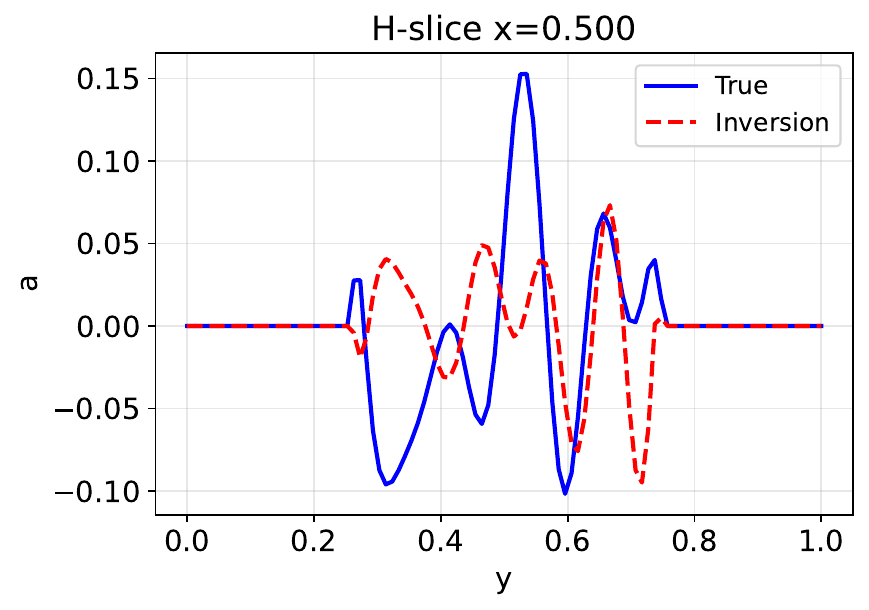} \\
{\scriptsize error=0.947}
\end{tabular} &
\begin{tabular}{c}
{\small $r=45^\circ, s=180^\circ$} \\
\includegraphics[height=2.4cm]{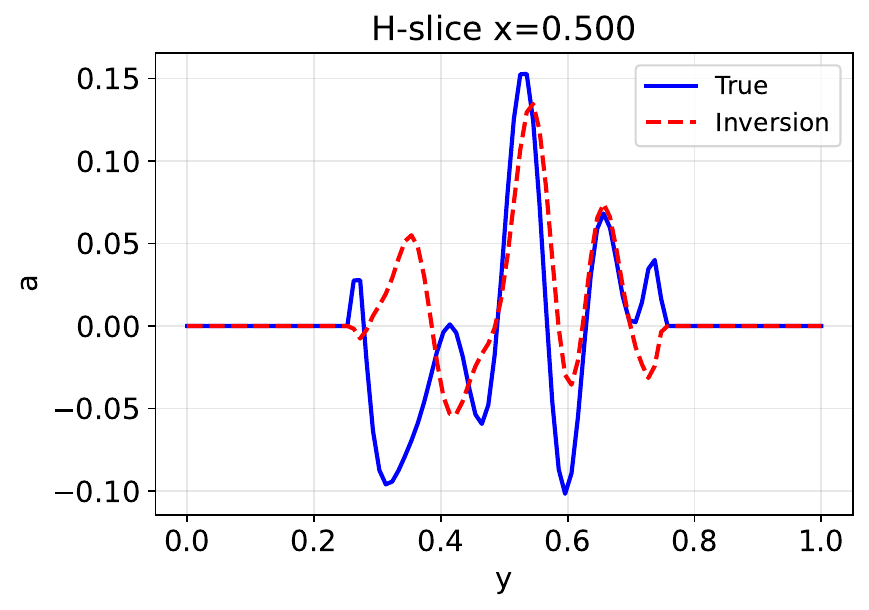} \\
{\scriptsize error=0.744}
\end{tabular} &
\begin{tabular}{c}
{\small $r=45^\circ, s=360^\circ$} \\
\includegraphics[height=2.4cm]{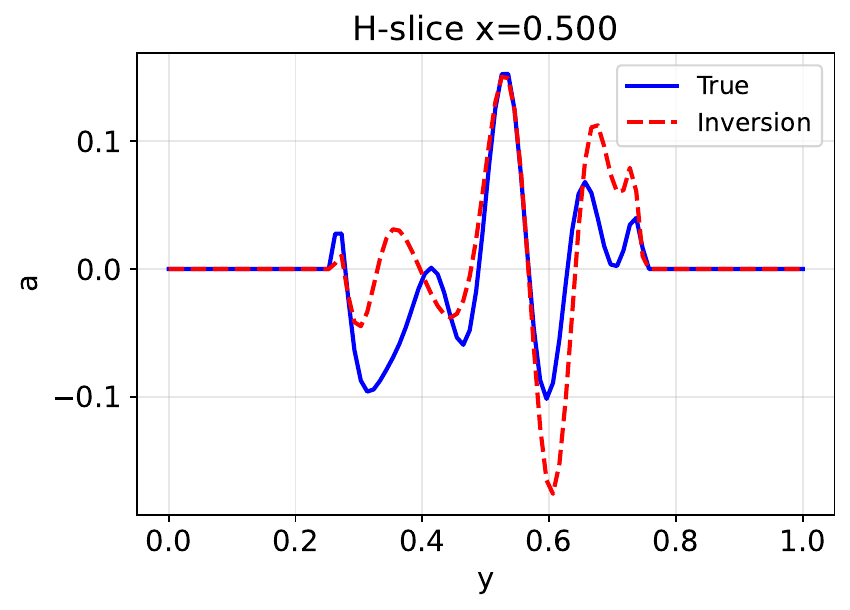} \\
{\scriptsize error=0.476}
\end{tabular}
\\

\begin{tabular}{c}
{\small $r=90^\circ, s=45^\circ$} \\
\includegraphics[height=2.4cm]{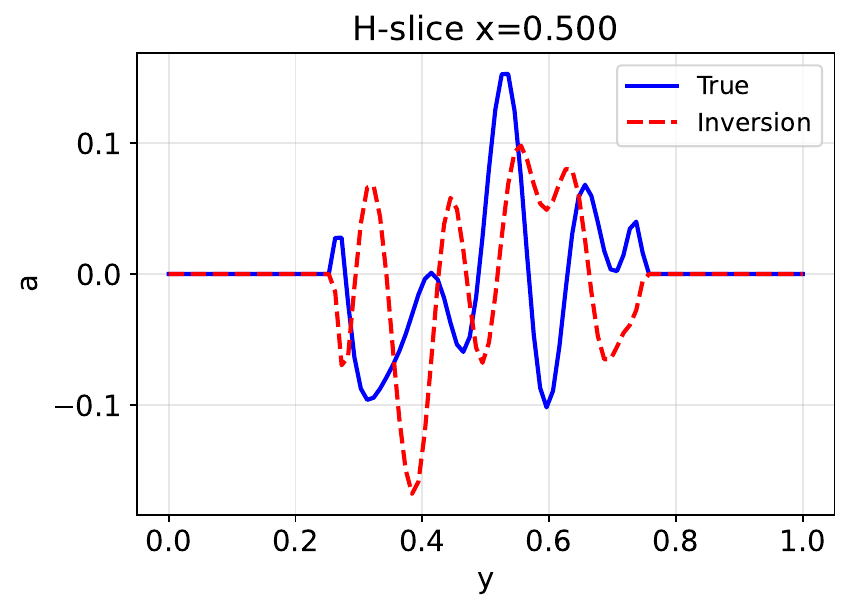} \\
{\scriptsize error=0.919}
\end{tabular} &
\begin{tabular}{c}
{\small $r=90^\circ, s=90^\circ$} \\
\includegraphics[height=2.4cm]{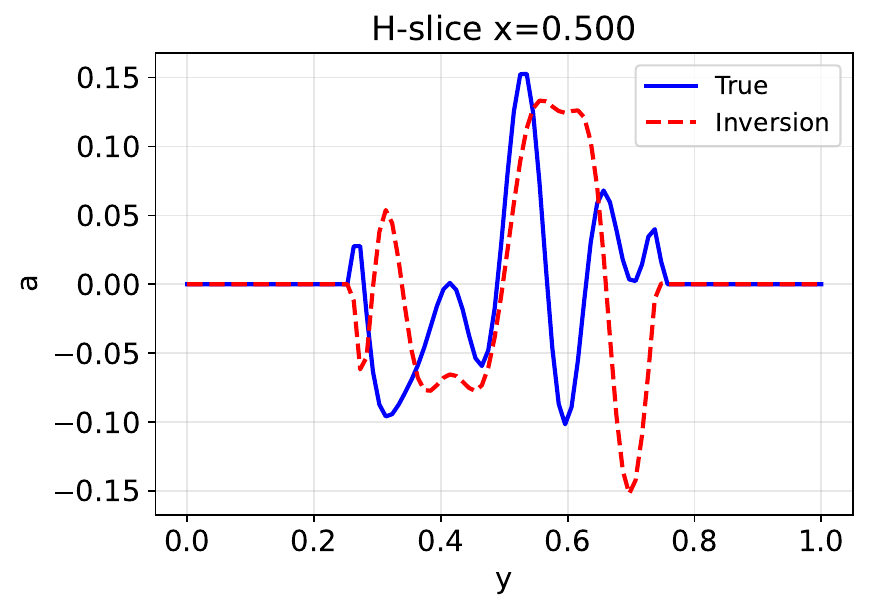} \\
{\scriptsize error=0.833}
\end{tabular} &
\begin{tabular}{c}
{\small $r=90^\circ, s=180^\circ$} \\
\includegraphics[height=2.4cm]{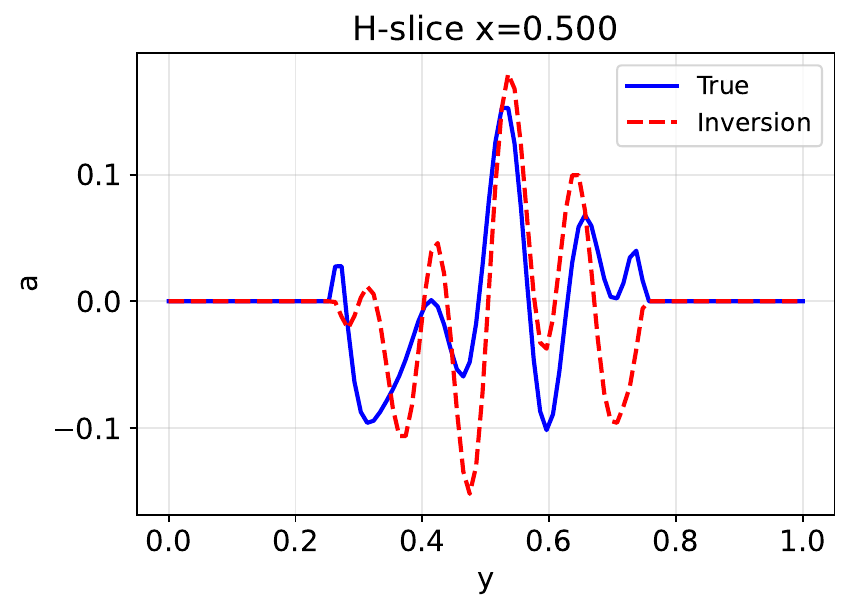} \\
{\scriptsize error=0.671}
\end{tabular} &
\begin{tabular}{c}
{\small $r=90^\circ, s=360^\circ$} \\
\includegraphics[height=2.4cm]{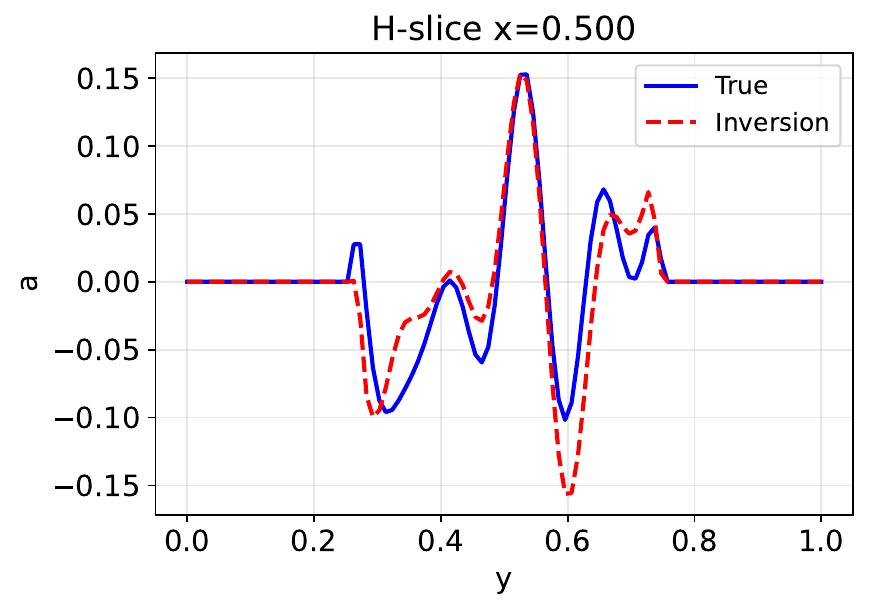} \\
{\scriptsize error=0.336}
\end{tabular}
\\

\begin{tabular}{c}
{\small $r=180^\circ, s=45^\circ$} \\
\includegraphics[height=2.4cm]{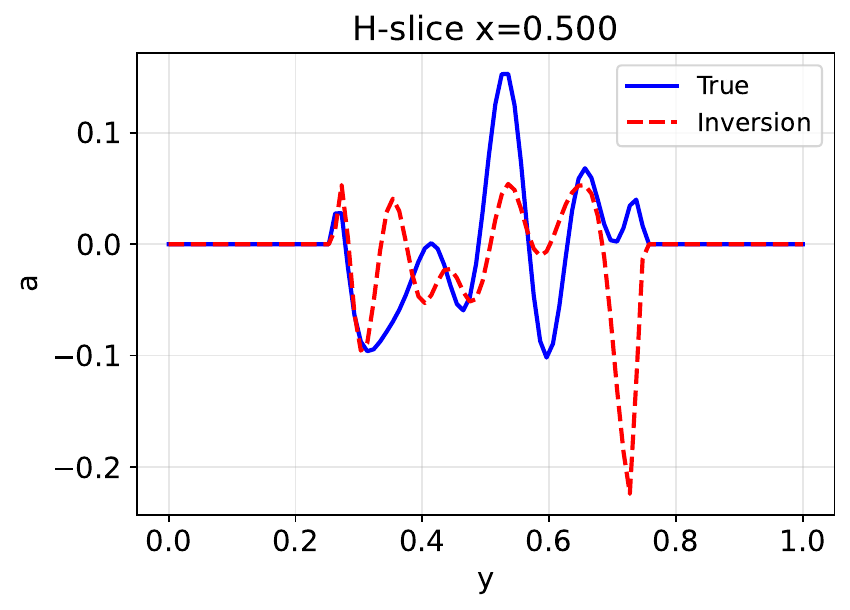} \\
{\scriptsize error=0.833}
\end{tabular} &
\begin{tabular}{c}
{\small $r=180^\circ, s=90^\circ$} \\
\includegraphics[height=2.4cm]{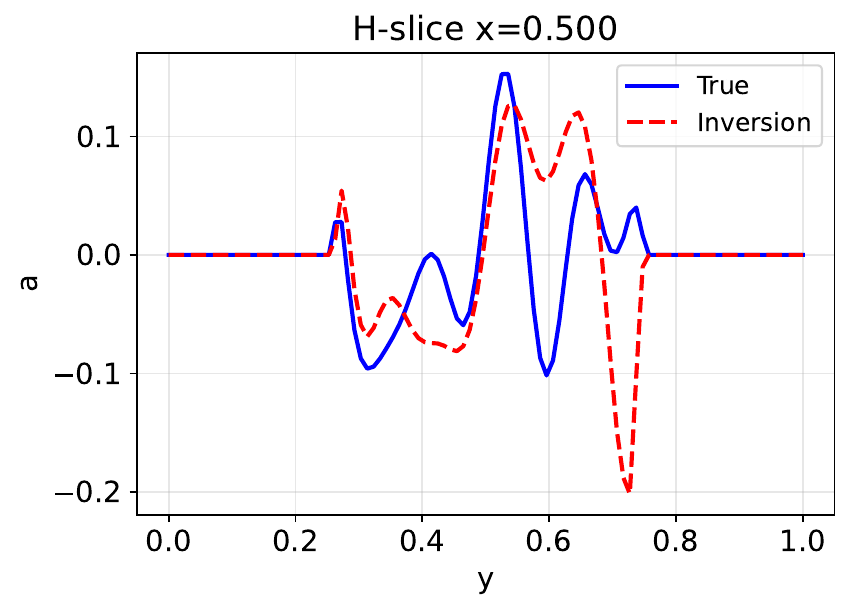} \\
{\scriptsize error=0.684}
\end{tabular} &
\begin{tabular}{c}
{\small $r=180^\circ, s=180^\circ$} \\
\includegraphics[height=2.4cm]{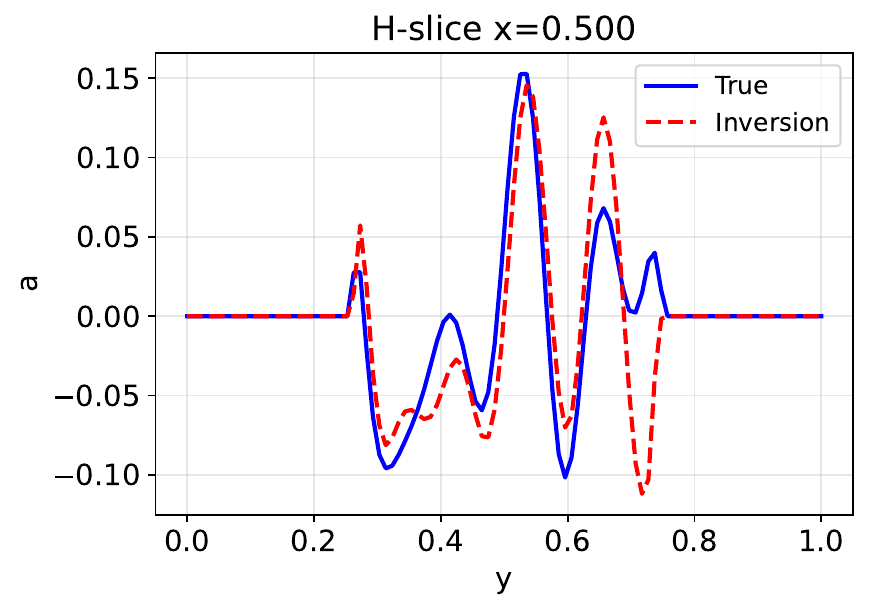} \\
{\scriptsize error=0.435}
\end{tabular} &
\begin{tabular}{c}
{\small $r=180^\circ, s=360^\circ$} \\
\includegraphics[height=2.4cm]{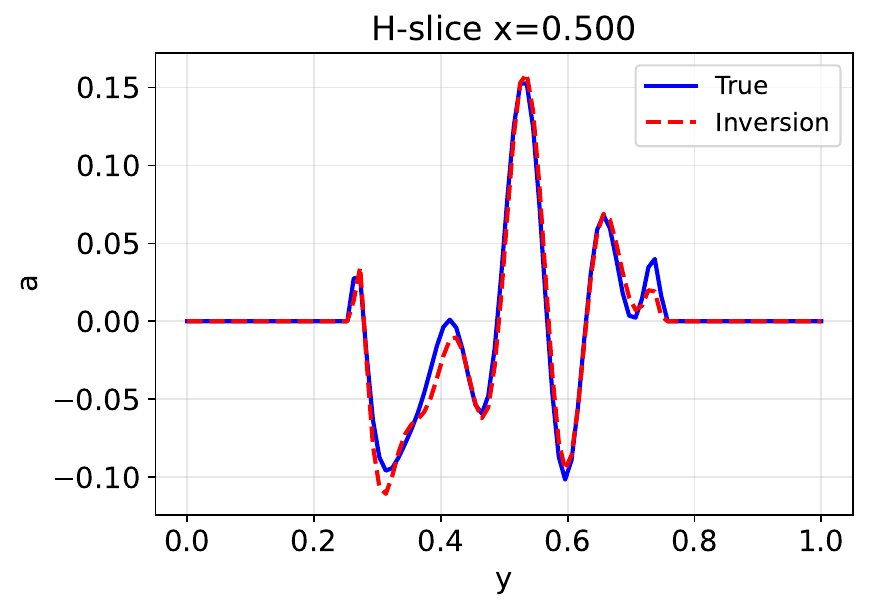} \\
{\scriptsize error=0.143}
\end{tabular}
\\

\begin{tabular}{c}
{\small $r=360^\circ, s=45^\circ$} \\
\includegraphics[height=2.4cm]{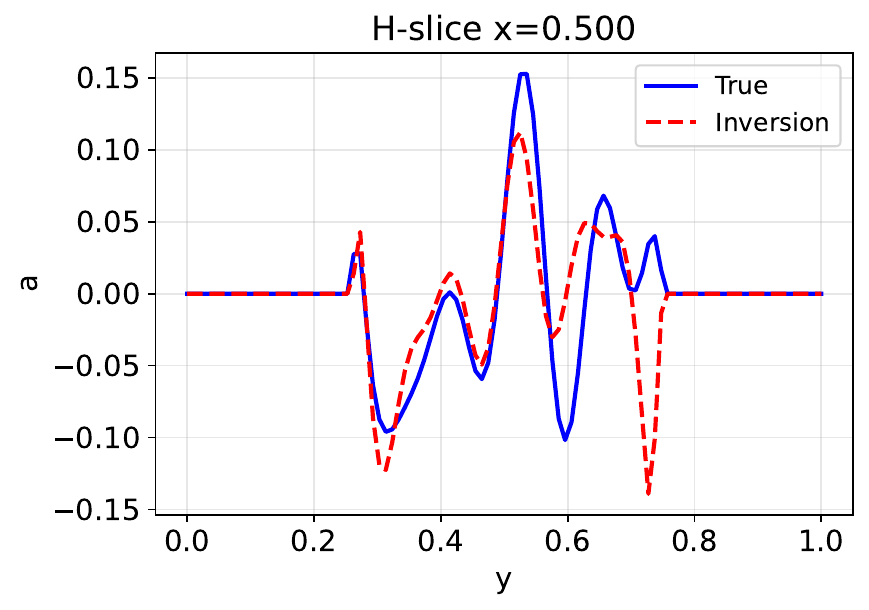} \\
{\scriptsize error=0.569}
\end{tabular} &
\begin{tabular}{c}
{\small $r=360^\circ, s=90^\circ$} \\
\includegraphics[height=2.4cm]{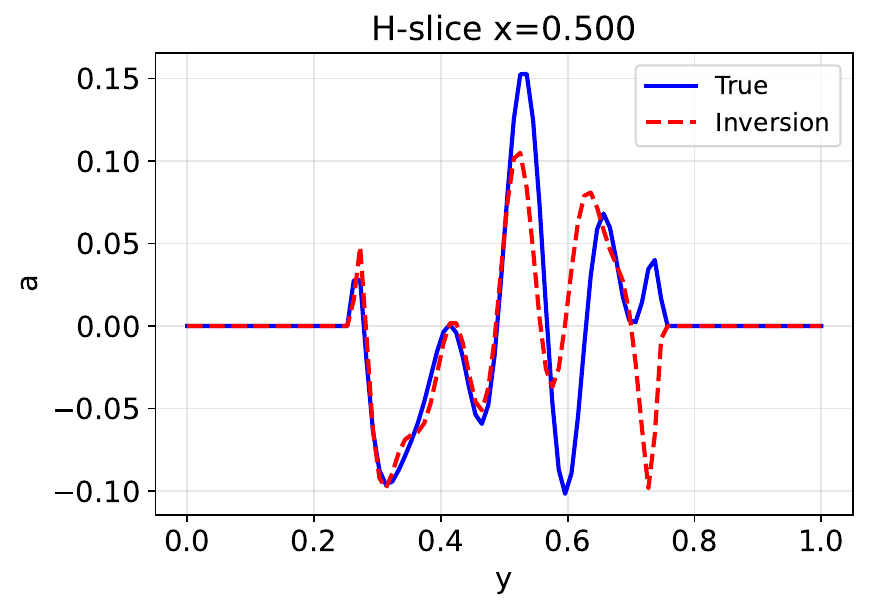} \\
{\scriptsize error=0.421}
\end{tabular} &
\begin{tabular}{c}
{\small $r=360^\circ, s=180^\circ$} \\
\includegraphics[height=2.4cm]{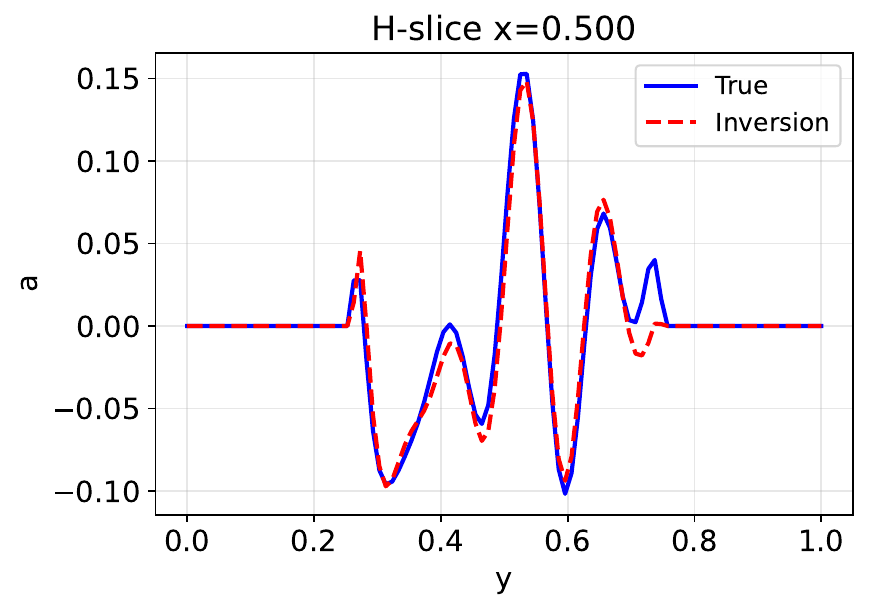} \\
{\scriptsize error=0.149}
\end{tabular} &
\begin{tabular}{c}
{\small $r=360^\circ, s=360^\circ$} \\
\includegraphics[height=2.4cm]{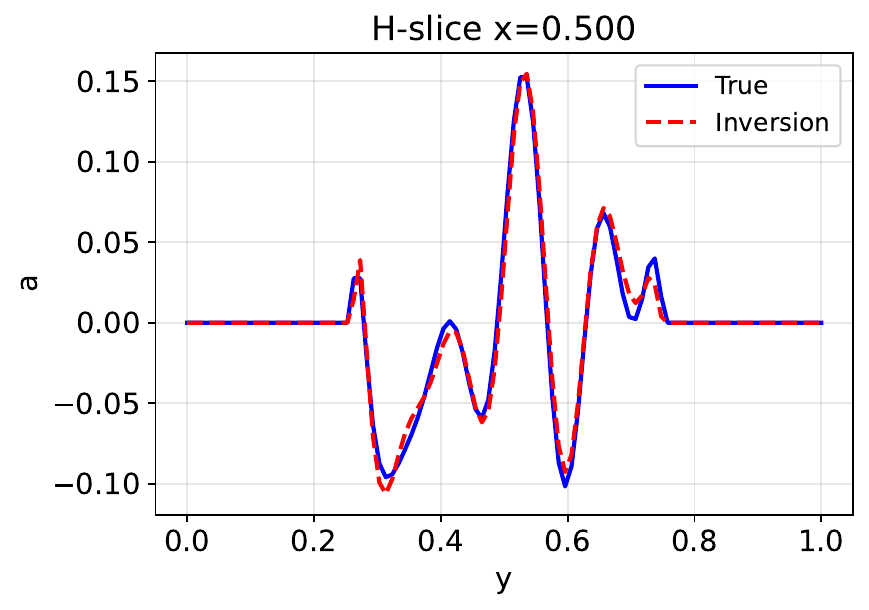} \\
{\scriptsize error=0.129}
\end{tabular}

\end{tabular}

\centering
\begin{tabular}{cc}

\begin{tabular}{c}
{\small Diagonal $r=45^\circ, s=360^\circ$} \\
\includegraphics[height=2.4cm]{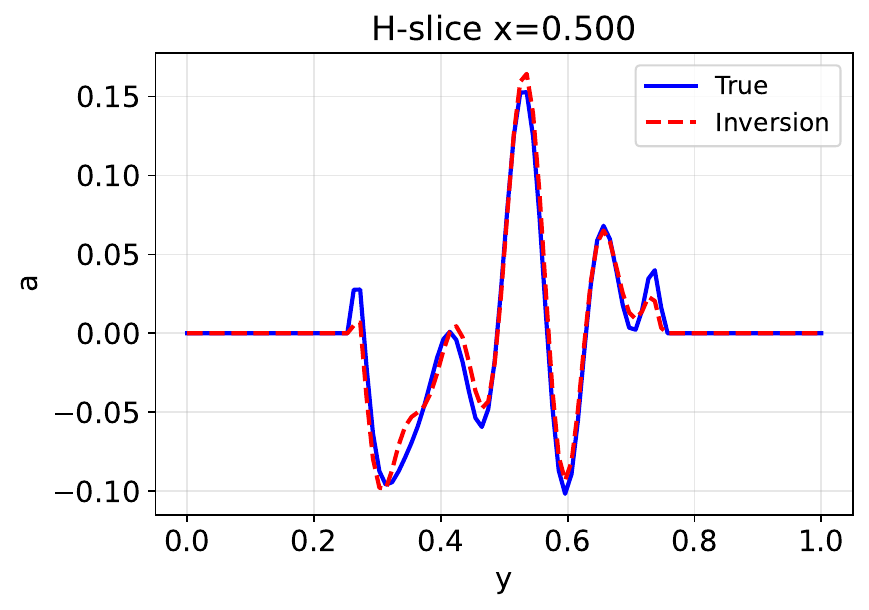} \\
{\scriptsize error=0.319}
\end{tabular}

&

\begin{tabular}{c}
{\small Diagonal $r=90^\circ, s=360^\circ$} \\
\includegraphics[height=2.4cm]{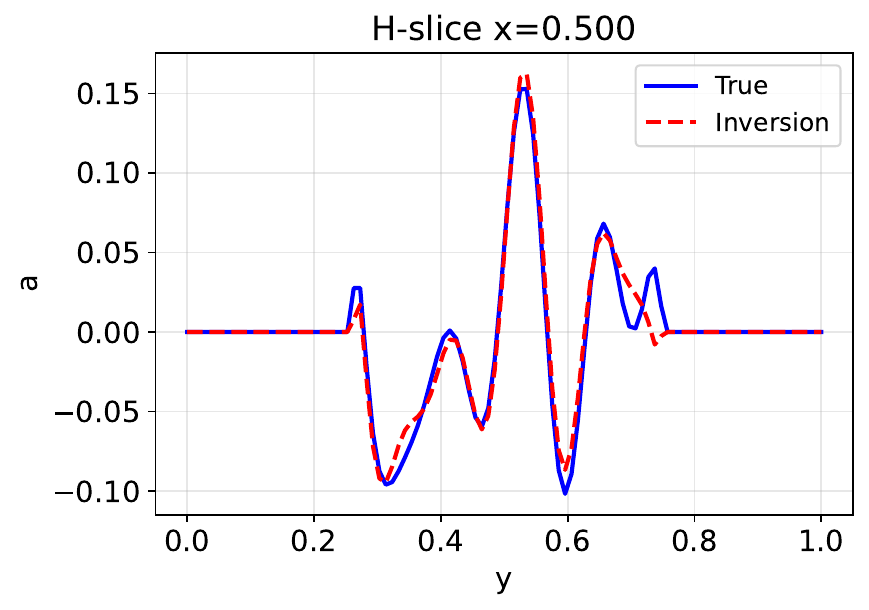} \\
{\scriptsize error=0.157}
\end{tabular}

\end{tabular}
\caption{
Reconstruction comparisons on the cross-section $x=0.5$ for different source and receiver apertures with error reported below each reconstruction.
}
\label{fig:reconstruction_hslices}
\end{figure}

Figure~\ref{fig:reconstruction_hslices} further shows cross-sectional comparisons along $x=0.5$. When both apertures are severely limited, the reconstructed profiles deviate significantly from the reference across the entire cross-section and exhibit strong directional bias. This behavior reflects the intrinsic ill-posedness induced by incomplete angular coverage, where parts of the domain are insufficiently illuminated or observed. As expected, widening the angular coverage of either sources or receivers consistently brings the reconstructed profile closer to the reference. In near full-aperture configurations, the reconstructed profiles closely match the reference solution, with only minor local discrepancies. These results underscore that sufficient spatial diversity in both illumination and observation is crucial for mitigating the intrinsic ill-posedness of limited-aperture inverse problems.

\noindent\textbf{Comparison of regularization strategies} To isolate the effect of each regularization strategy, all reconstructions in this subsection use clean (noise-free) observations at wavenumber $k=60$ with the forward model taken as the MscaleFNO surrogate trained in Section~\ref{sec:forward}. We compare three strategies: (1) the proposed Plug-and-Play (PnP) approach, where the EDM denoiser replaces the proximal operator within the iterative reconstruction, (2) standard Tikhonov regularization, which imposes an explicit $\ell_2$-smoothness penalty, and (3) unregularized inversion. In the PnP framework, no explicit regularization functional is required; instead, the learned denoiser implicitly constrains the solution to remain near the data manifold on which the forward surrogate was trained.

\begin{figure}[htbp]
\centering

\begin{subfigure}{0.9\textwidth}
    \centering
    \includegraphics[width=0.4\linewidth]{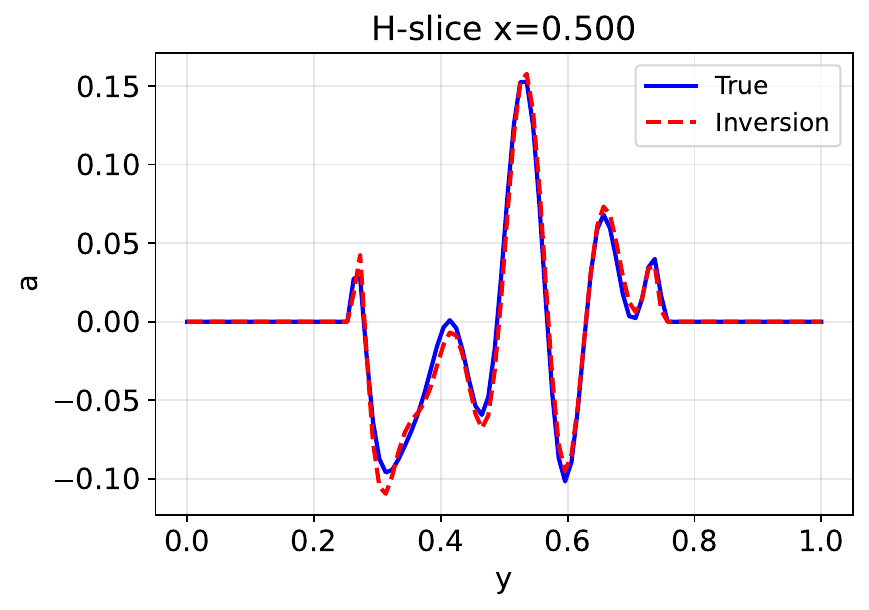}
    \includegraphics[width=0.4\linewidth]{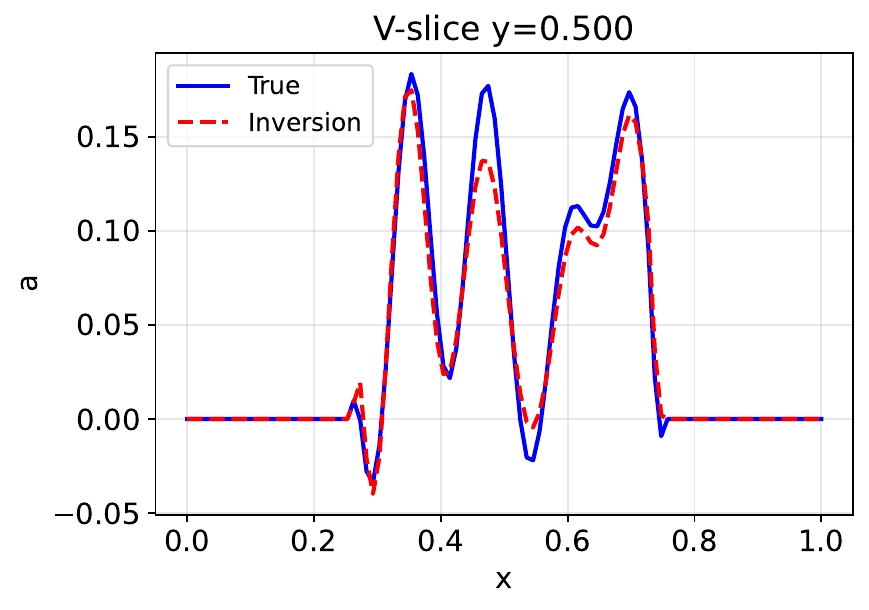}
    \caption{PnP with EDM denoiser: relative error = 0.122}
\end{subfigure}

\vspace{0.5em}

\begin{subfigure}{0.9\textwidth}
    \centering
    \includegraphics[width=0.4\linewidth]{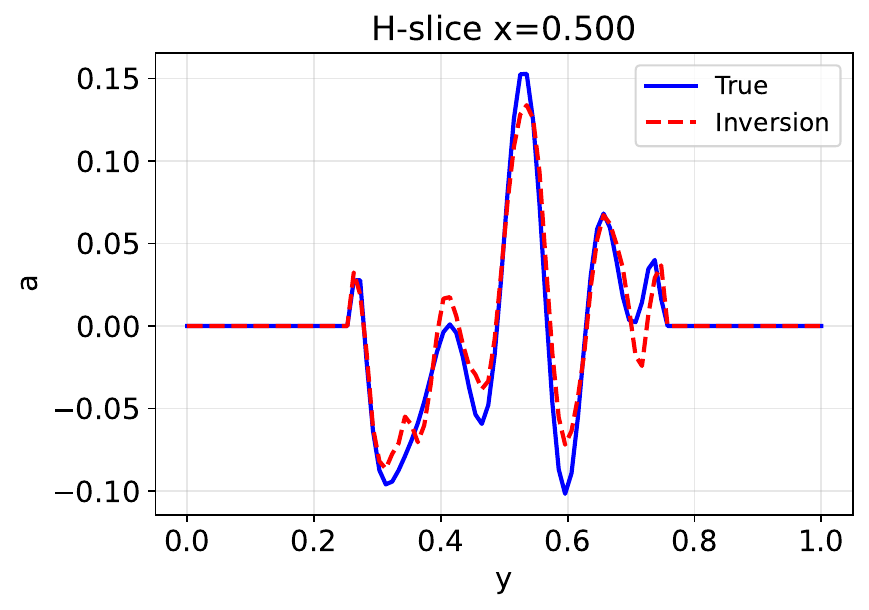}
    \includegraphics[width=0.4\linewidth]{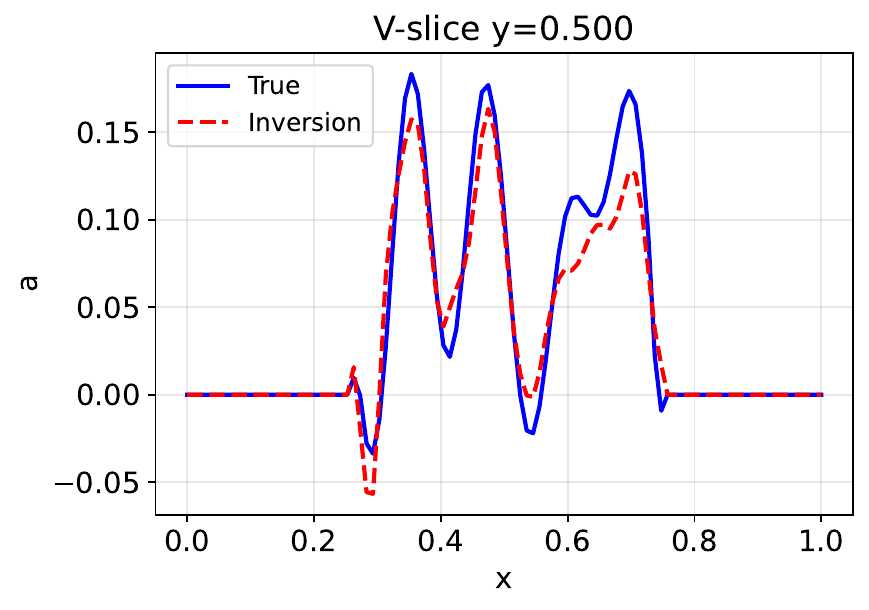}
    \caption{Tikhonov regularization: relative error = 0.282}
\end{subfigure}

\vspace{0.5em}

\begin{subfigure}{0.9\textwidth}
    \centering
    \includegraphics[width=0.4\linewidth]{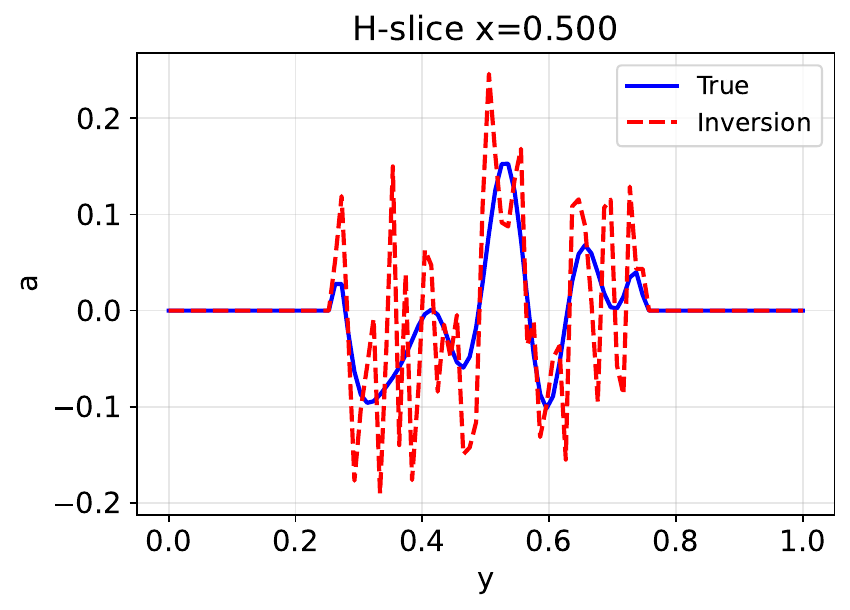}
    \includegraphics[width=0.4\linewidth]{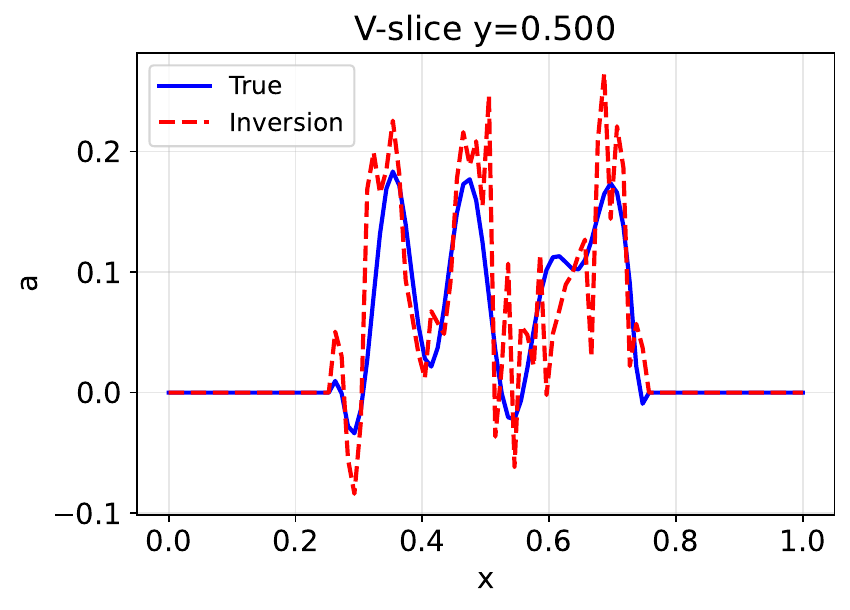}
    \caption{No regularization: relative error = 0.549}
\end{subfigure}

\caption{Reconstruction results under different strategies at $k=60$ using clean observations. For each method, the left panel shows the horizontal slice at $y=0.5$ and the right panel shows the vertical slice at $x=0.5$. The PnP reconstruction preserves oscillatory features with minimal artifacts, Tikhonov regularization oversmoothes sharp transitions, and the unregularized solution is dominated by spurious high-frequency oscillations.}
\label{fig:reg_comparison}
\end{figure}

As shown in Figure~\ref{fig:reg_comparison}, the unregularized inversion amplifies high-frequency errors inherent in the forward surrogate, producing spurious oscillations that mask the true medium structure. Tikhonov regularization suppresses these artifacts through explicit smoothness constraints, but damps genuine oscillatory features in the process, with peaks and troughs visibly attenuated. In contrast, the PnP approach with the EDM denoiser effectively removes spurious oscillations while faithfully preserving the high-frequency components that define the medium. This follows directly from how the PnP framework operates. Instead of imposing a generic smoothness penalty, the denoiser projects each iterate onto a learned manifold of physically admissible fields, which allows it to separate genuine structure from artifact.

\begin{figure}[htbp]
\centering
\includegraphics[width=0.65\linewidth]{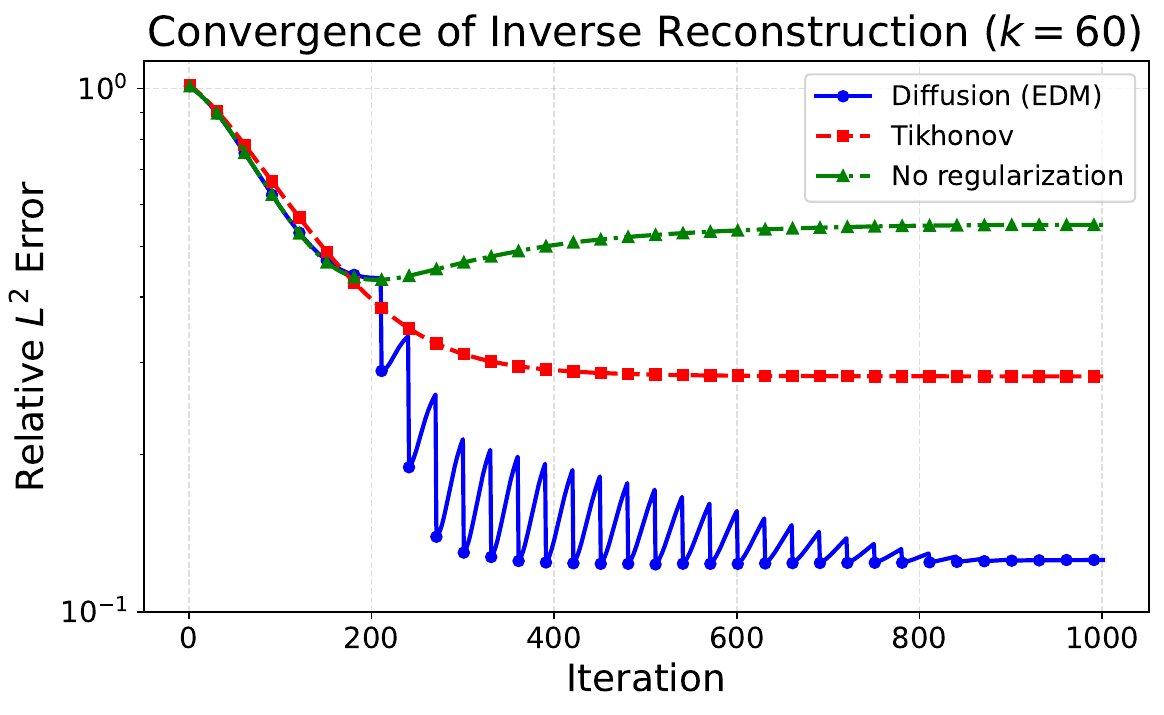}
\caption{Convergence behavior of different strategies at $k=60$ under clean observations. The relative $L^2$ error is computed with respect to the ground truth.}
\label{fig:reg_convergence}
\end{figure}

The convergence profiles in Figure~\ref{fig:reg_convergence} reinforce these observations. The unregularized method descends briefly before stagnating at a high error, a consequence of the problem's intrinsic ill-posedness. Tikhonov regularization converges stably but saturates at a moderate error, constrained by its preference for smooth solutions.

The PnP method converges in a sawtooth pattern: each sharp drop coincides with a denoising step that pulls the estimate back onto the learned manifold of admissible fields, while the intervening data-fidelity updates may temporarily nudge the error upward. This oscillatory dynamics helps the solver escape poor local minima that trap the unregularized and Tikhonov schemes, and eventually drives the error well below the other two.

These results point to a fundamental mismatch between explicit smoothness priors and high-wavenumber inverse problems. The signal of interest lives in the very oscillatory components that such regularizers penalize. The PnP framework sidesteps this conflict by replacing the hand-crafted proximal operator with a data-driven denoiser, yielding reconstructions that are both stable and faithful to the fine-scale structure of the medium.

\subsection{Ablation study}
\subsubsection{Spectral analysis of multiscale branches}

To better understand the mechanism of the proposed MscaleFNO, we analyze the spectral characteristics of different branches in the network.
Specifically, we consider the operator learned by MscaleFNO with $k=60$. The model consists of four scale branches, with initial scales $\bm{c} = \{1, 100, 500, 2000\}$. Both the weights and the scales are trainable. After training, the learned scales become $\bm{c} = \{2.99, 91.80, 491.73, 1991.79\}$, showing that the scales remain largely stable during optimization.

\begin{figure}[t!]
\centering
\begin{subfigure}{0.48\textwidth}
    \centering
    \includegraphics[width=\linewidth]{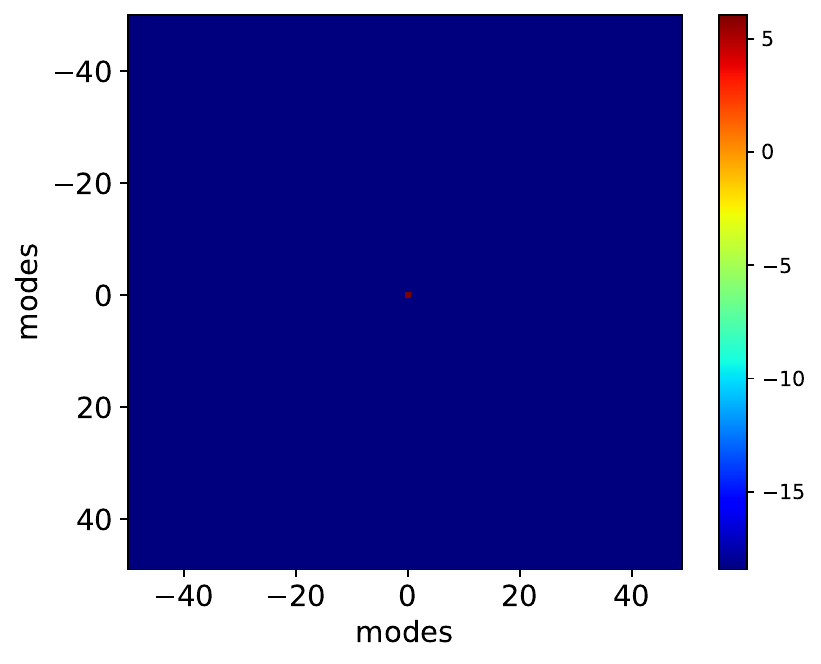}
    \caption{scale = $2.99$}
\end{subfigure}
\begin{subfigure}{0.48\textwidth}
    \centering
    \includegraphics[width=\linewidth]{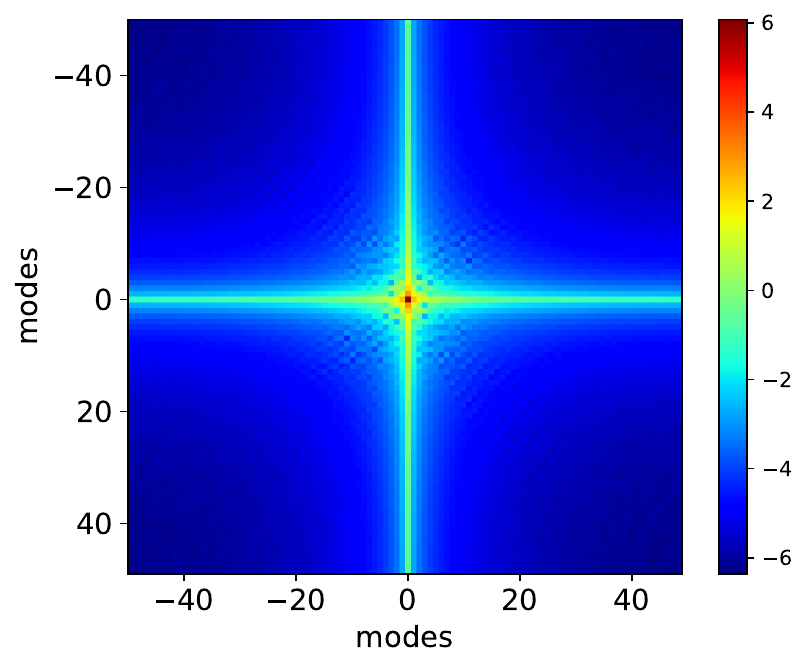}
    \caption{scale = $91.80$}
\end{subfigure}
\\
\begin{subfigure}{0.48\textwidth}
    \centering
    \includegraphics[width=\linewidth]{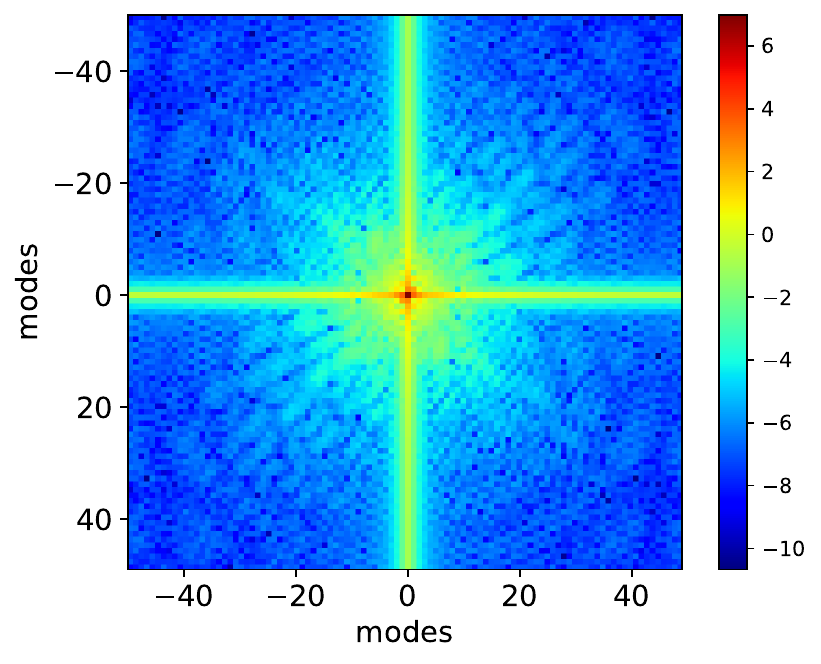}
    \caption{scale = $491.73$}
\end{subfigure}
\begin{subfigure}{0.48\textwidth}
    \centering
    \includegraphics[width=\linewidth]{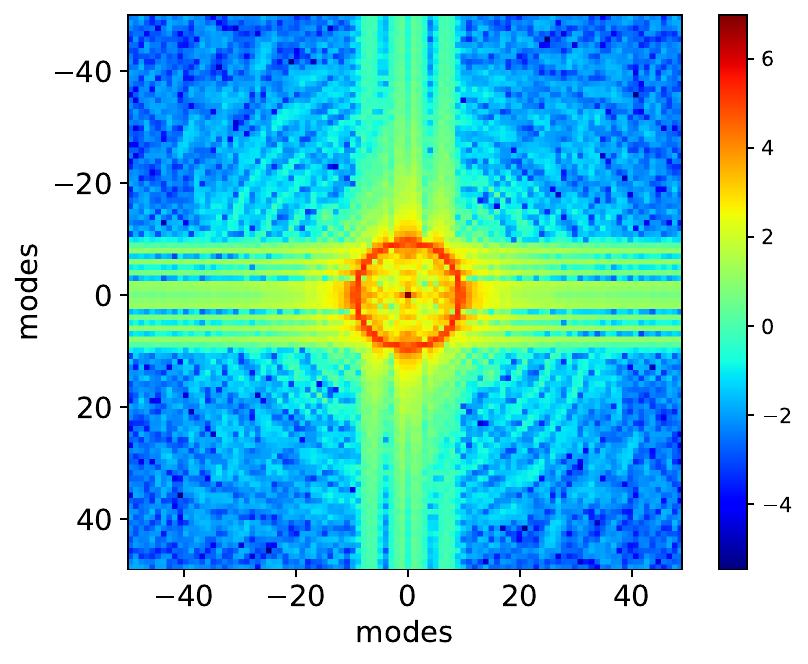}
    \caption{scale = $1991.79$}
\end{subfigure}

\caption{Fourier spectra of the outputs from different branches in MscaleFNO (trained with $k=60$). Each subfigure corresponds to one branch with its learned scale. Distinct spectral patterns are observed across branches.}
\label{fig:branch_spectrum}
\end{figure}

We then examine the Fourier spectrum of the outputs from each scale branch. Although each branch operates at a different spatial scale, we directly visualize their frequency spectra to investigate how their outputs differ in the frequency domain.
Figure~\ref{fig:branch_spectrum} shows the magnitude of the Fourier coefficients for different branches. A clear trend can be observed across scales: the branch with the smallest scale ($c=2.99$) exhibits relatively weak spectral responses, while branches with larger scales display increasingly rich and structured frequency patterns. As the scale increases, the spectral energy becomes more pronounced and spreads over a wider range of modes.
Such behavior aligns with the design motivation of MscaleFNO and suggests that the model is not merely increasing capacity, but learning a structured representation across frequencies. Moreover, this frequency-wise differentiation helps alleviate the optimization difficulties observed in standard FNO, where different frequency components are entangled within a single representation.

Overall, these results indicate that the proposed multiscale design promotes a degree of frequency separation, which is beneficial for modeling high-frequency phenomena.

\subsubsection{Effect of Fourier modes and network width in FNO}

To ensure a fair comparison and rule out the possibility that the poor performance of FNO at high frequencies is caused by insufficient spectral resolution, we conduct an ablation study by increasing the number of Fourier modes to the maximum allowable value. Specifically, given the training resolution of $101 \times 101$, we set the number of Fourier modes to 50, which corresponds to (or closely approaches) the Nyquist frequency limit of the discretization. This ensures that the model has access to the highest resolvable frequency components, eliminating spectral truncation as a limiting factor.
In addition, we vary the network width, another key hyperparameter in FNO, and evaluate performance across different configurations. The tested widths are 16, 32, 48, and 64. For each configuration, we evaluate the learned forward operator across a range of wavenumbers $k = 10, 20, 30, 40, 50, 60$.

\begin{table}[H]
\centering
\caption{Relative error of FNO with maximum Fourier modes (50) under different network widths and wavenumbers $k$.}
\label{tab:fno_ablation}
\begin{tabular}{c|cccccc}
\hline
\textbf{Width \& $k$} & \textbf{10} & \textbf{20} & \textbf{30} & \textbf{40} & \textbf{50} & \textbf{60} \\
\hline
16 & 0.0242 & 0.0377 & 0.0554 & 0.9999 & 0.9999 & 0.9999 \\
32 & 0.0185 & 0.0210 & 0.0384 & 0.9999 & 0.9999 & 0.9999 \\
48 & 0.0136 & 0.0291 & 0.0415 & 0.1355 & 0.9999 & 0.9999 \\
64 & 0.0159 & 0.0314 & 0.0381 & 0.9999 & 0.9999 & 0.9999 \\
\hline
\end{tabular}
\end{table}

From Table~\ref{tab:fno_ablation}, we observe a clear transition in the behavior of FNO as the wavenumber $k$ increases. For relatively small values of $k$ (e.g., $k \leq 30$), FNO achieves low relative errors across all tested widths, indicating that the model is capable of learning operators with moderate frequency content and nonlinearity. However, as $k$ increases beyond this regime, the performance deteriorates significantly. In particular, for $k \geq 40$, FNO starts to fail to learn the operator, with the relative error saturating at a high level (close to $1.0$ in most cases), suggesting that the training process stagnates.
To investigate whether this failure is due to insufficient model capacity, we increase the network width up to 64 while keeping the number of Fourier modes at the Nyquist limit. We observe that increasing the width from 16 to 48 leads to partial improvement (e.g., error $0.1355$ at width 48 for $k=40$), indicating that optimization can be temporarily facilitated with a larger model. However, further increasing the width does not lead to additional gains; instead, the training again fails to converge, and the error returns to a saturated regime.
These results indicate that the failure of FNO in high-frequency regimes is not solely due to insufficient spectral resolution or limited model capacity. Even with maximal Fourier modes and large network width, FNO remains difficult to optimize when the target operator exhibits strong nonlinearity and high-frequency oscillations. This suggests that the limitation is not primarily in approximation capability, but is likely related to the difficulty of optimizing representations involving strongly coupled multi-frequency components.

In contrast, under the same training settings (e.g., learning rate, batch size, optimizer, and scheduler), MscaleFNO achieves this without ever retaining the full spectrum on a single, high-resolution grid. This indicates that the limitation of FNO is intrinsically related to its single-scale representation, which entangles different frequency components and leads to optimization instability. The multiscale design of MscaleFNO, on the other hand, provides a more structured representation that better aligns with the underlying frequency characteristics. By distributing the learning across multiple scales, the model mitigates the coupling between different frequency components and improves optimization behavior. As a result, it not only improves accuracy but also significantly enhances training robustness and efficiency under limited model capacity.

\section{Conclusion}

In this work, we have proposed a full-waveform inversion framework 
for oscillatory inverse medium problems governed by the Helmholtz 
equation, centered on MscaleFNO as the core differentiable surrogate. 
Its multiscale spectral architecture directly addresses the spectral 
bias of standard FNO, enabling accurate resolution of rapidly 
oscillatory wavefields across high-wavenumber regimes and providing 
a faithful differentiable surrogate for the Helmholtz forward map. 
This high-frequency fidelity is the key enabler of reliable 
gradient-based inversion for fine-scale medium recovery, which 
standard FNO fails to achieve at large wavenumbers.

Numerical experiments confirm that MscaleFNO achieves substantially 
more accurate forward predictions than standard FNO in high-wavenumber 
regimes, and that this advantage directly translates to improved 
inversion quality for high-frequency media. An EDM-based plug-and-play 
prior is further incorporated as an auxiliary regularization to 
improve robustness under noisy and limited-aperture conditions, 
where severely restricted apertures limit fine-scale recovery but 
main oscillatory structures are reconstructed more reliably as the 
observation aperture increases.






\bibliographystyle{plain}

\end{document}